\documentclass[reqno,11pt]{amsart}

\usepackage{amsthm, amsmath, amssymb, stmaryrd}
\usepackage[utf8]{inputenc}
\usepackage[T1]{fontenc}

\usepackage[hyphenbreaks]{breakurl}
\usepackage[hyphens]{url}
\usepackage{systeme}
\usepackage[shortlabels]{enumitem}
\usepackage[hidelinks]{hyperref}
\usepackage{microtype}

\usepackage{bm}
\usepackage[margin=1in]{geometry}

\usepackage[textsize=scriptsize,backgroundcolor=orange!5]{todonotes}

\usepackage[noabbrev,capitalize,sort]{cleveref}
\crefname{equation}{}{}
\crefname{enumi}{}{}
\numberwithin{equation}{section}

\usepackage{mathtools}
\newtheorem{theorem}{Theorem}[section]
\newtheorem{proposition}[theorem]{Proposition}
\newtheorem{lemma}[theorem]{Lemma}
\newtheorem{claim}[theorem]{Claim}
\newtheorem{corollary}[theorem]{Corollary}

\theoremstyle{definition}
\newtheorem{definition}[theorem]{Definition}

\newtheorem{example}[theorem]{Example}

\theoremstyle{remark}
\newtheorem*{remark}{Remark}

\newcommand{\abs}[1]{\left\lvert#1\right\rvert}

\newcommand{\norm}[1]{\left\lVert#1\right\rVert}

\DeclareMathOperator{\supp}{supp}   
\DeclareMathOperator{\ex}{ex}

\newcommand*{\eqdef}{\stackrel{\mbox{\normalfont\tiny def}}{=}}

\newcommand{\FF}{\mathbb{F}}

\newcommand{\RR}{\mathbb{R}}
\newcommand{\NN}{\mathbb{N}}

\newcommand*{\PP}{\mathbb{P}}

\newcommand{\cF}{\mathcal F}

\newcommand{\cE}{\mathcal E}

\newcommand{\bi}{\mathbf{i}}
\newcommand{\HH}{\mathbb{H}}

\newlength{\hght}

\makeatletter
\newcommand\thankssymb[1]{\textsuperscript{\@fnsymbol{#1}}}
\makeatother

\author[Ting-Wei Chao]{Ting-Wei Chao\thankssymb{1}}
\author[Hung-Hsun Hans Yu]{Hung-Hsun Hans Yu\thankssymb{2}}

\thanks{\thankssymb{1}Department of Mathematics, Massachusetts Institute of Technology, Cambridge, MA, USA. Email: {\tt twchao@mit.edu}}

\thanks{\thankssymb{2}Department of Mathematics, Princeton University, Princeton, NJ 08544\@.  Email: {\tt hansonyu@princeton.edu}}


\title{When entropy meets Tur\'an: \linebreak new proofs and hypergraph Tur\'an results}

\begin{document}

\maketitle

\begin{abstract}
In this paper, we provide a new proof of a density version of Tur\'an's theorem.
We also rephrase both the theorem and the proof using entropy.
With the entropic formulation, we show that some naturally defined entropic quantity is closely connected to other common quantities such as Lagrangian and spectral radius.
In addition, we also determine the Tur\'an density for a new family of hypergraphs, which we call tents.
Our result can be seen as a new generalization of Mubayi's result on the extended cliques.
\end{abstract}

\section{Introduction}
For any $k$-graph (i.e.\ $k$-uniform hypergraph) $F$, its \emph{Tur\'an number} $\ex(n,F)$ is the maximum number of edges in an $F$-free $k$-graph $G$ on $n$ vertices.
Here, $G$ is $F$-free if it contains no subgraph (not necessarily induced) isomorphic to $F$.
The study of Tur\'an numbers was initiated by Tur\'an \cite{Turan41}, who first considered the case where $k=2$ and $F$ is the complete graph $K_{r+1}$ on $(r+1)$ vertices.
There, Tur\'an showed that $\ex(n,F)$ is maximized by the balanced complete $r$-partite graph $T_{n,r}$, which we now refer to as the Tur\'an graph.
Tur\'an's foundational work has motivated subsequent works on related problems, driving continuing research in extremal graph theory.

The general Tur\'an problem is fairly understood when $k=2$.
Although the exact value of $\ex(n,F)$ is not known for general graphs $F$, the celebrated Erd\H{o}s--Stone theorem \cite{ES46} asserts that $\ex(n,F) = \left(1-\frac{1}{r}+o(1)\right)\binom{n}{2}$ if $\chi(F) = r+1$, where $T_{n,r}$ is an asymptotic extremizer.
If we define the \emph{Tur\'an density} to be
\[\pi(F) = \lim_{n\to\infty}\frac{\ex(n,F)}{\binom{n}{k}}\]
for a $k$-graph $F$, then the Erd\H{o}s--Stone theorem can be rephrased as $\pi(F) = 1-\frac{1}{\chi(F)-1}$ when $F$ is a graph.
It is worth pointing out that when $\chi(F)=2$, Erd\H{o}s--Stone gives that $\pi(F)=0$, showing that $\ex(n,F)$ is subquadratic but does not determine the asymptotic behavior of $\ex(n,F)$.
Despite lots of effort, there are still many interesting open problems regarding the asymptotic behavior of $\ex(n,F)$ when $F$ is bipartite.
However, in this paper, we will focus on the non-degenerate case where $\pi(F)>0$.

Given how much we know about Tur\'an numbers and Tur\'an densities of graphs, it might be surprising how little we know about hypergraph Tur\'an problems.
In fact, the exact value of $\pi(F)$ is still unknown even for $F=K_4^{(3)}$, the $3$-uniform clique on $4$ vertices.
Tur\'an \cite{Tur90} showed that $\pi(K_4^{(3)})\geq \frac{5}{9}$ and conjectured that it is actually an equality.
However, proving this conjecture still seems hard to date, and the current best upper bound $\pi(F)\leq 0.561666$ was obtained by Razborov \cite{Raz10} using flag-algebraic computation, which was later verified by \cite{BT11} and \cite{F-RV13}.
The difficulty comes from the fact that hypergraph Tur\'an problems have drastically different behaviors from the graph case.
For example, there is a large family of constructions all showing $\pi(K_4^{(3)})\geq \frac{5}{9}$ given in \cite{Kos82} (also see \cite{F-D-F88}).
In comparison, the Erd\H{o}s--Simonovits theorem states that any asymptotic extremizer of $\pi(K_{r+1})$ should be close to $T_{n,r}$.
We will discuss other interesting phenomena for hypergraph Tur\'an problems in \cref{subsec:hypergraph-turan-density}.

The aim of this paper is to find inspiration for new ways to approach hypergraph Tur\'an problems by examining our new proof of the density Tur\'an theorem.
\begin{theorem}[Density Tur\'an theorem]\label{thm:density-turan}
    For any positive integer $r$, 
    \[\pi(K_{r+1})=1-\frac{1}{r}.\]
\end{theorem}
Our new proof leads to new hypergraph Tur\'an results regarding hypergraphs that we call ``tents'', which generalize Mubayi's result \cite{Mub06} on the extended cliques.
We will introduce our results and related work in more detail in \cref{subsec:hypergraph-turan-density}.

Before diving into hypergraph Tur\'an problems, we will first give a quick overview of known proofs of Tur\'an's theorem.
We will then introduce the entropy method, which we use to rephrase both the theorem statement and our proof.
Then we will mention our hypergraph Tur\'an results that can be obtained using the new perspective, which can be thought of as one of our main results.

\subsection{Proofs of Tur\'an's theorem}
Tur\'an's original proof \cite{Turan41} works by a clever induction on the number of vertices by removing a $K_r$ from the graph.
Erd\H{o}s \cite{Erdos70} later provided another proof that modified the graph step by step, maintaining the $K_{r+1}$-freeness and making the graph complete multipartite at the end.
This method has the benefit that it is easier to see that the Tur\'an graph $T_{n,r}$ is the extremizer.
A proof of the same spirit is a folklore proof that proceeds with symmetrization (also known now as Zykov Symmetrization as this trick was used by Zykov \cite{Zyk49,Zyk52} in his work).
The proof modifies the graph by taking two non-adjacent vertices, and replacing one with another (see \cite[Chapter 41]{AZ18}).
Unfortunately, all those proofs do not easily generalize to hypergraphs as they all use properties of graphs crucially.

One proof that looks entirely different from the previous proofs is by applying the Caro--Wei theorem, which is due to Alon and Spencer \cite{AS00}.
The Caro--Wei theorem, independently proven by Caro \cite{Caro79} and Wei \cite{Wei81}, gives a lower bound on the independence number of a graph $G$ based on its degree sequence.
The standard proof of the Caro--Wei theorem is a nice probabilistic argument, which can be found in \cite{AS00}.
By taking the complement and an application of Cauchy--Schwarz, the density Tur\'an theorem immediately follows from Caro--Wei.
However, this argument does not generalize well to higher uniformities---although the Caro--Wei theorem can be extended to hypergraphs (see \cite{CT91}), applying the inequality on the complement no longer gives tight hypergraph Tur\'an results.

Another proof that is seemingly different from all the above is a proof due to Motzkin and Straus \cite{MS65}.
Their proof relies crucially on a quantity called \emph{Lagrangian}.
The Lagrangian $L(G)$ of a graph $G=(V,E)$ is defined as
\[\max \sum_{\{u,v\}\in E}x_ux_v \textup{ subj. to } x_v\geq 0\quad\forall v\in V\textup{ and }\sum_{v\in V}x_v=1.\]
Despite its somewhat long definition, it is a natural quantity to consider in the context of Tur\'an problems.
To see this, let $N$ be some large positive integer.
Consider the \emph{blowup} of $G$ obtained by putting in $(x_v+o(1))N$ copies of each vertex $v\in V$ so that there are $N$ vertices in total, where $(x_v)_{v\in V}$ is the extremizer for the Lagrangian.
Then there are $(L(G)+o(1))N^{2}$ edges in the blowup.
On the other hand, it is clear that $\abs{E}\leq L(G)\abs{V}^2$, which shows that the density Tur\'an theorem is equivalent to that $L(G)\leq \frac{1}{2}\left(1-\frac{1}{r}\right)$ for every $K_{r+1}$-free graph $G$.
Motzkin and Straus' idea is that if $u$ and $v$ are not adjacent, then there is an extremizer with either $x_u=0$ or $x_v=0$ for $L(G)$.
Therefore if $G$ is $K_{r+1}$-free, then there is an extremizer with support of size at most $r$.
A simple application of Cauchy--Schwarz then concludes the proof.
Despite its algebraic look, this proof is actually similar to Zykov Symmetrization in spirit.

It is natural to generalize graph Lagrangian to hypergraph Lagrangian.
For any $k$-graph $G=(V,E)$, its \emph{hypergraph Lagrangian} $L(G)$ is defined as the maximum of $\sum_{\{u_1,\ldots, u_k\}\in E}x_{u_1}\cdots x_{u_v}$ under the same conditions.
As before, when each $v\in V$ is blown-up to $(x_v+o(1))N$ vertices where $(x_v)_{v\in V}$ is the extremizer for the Lagrangian, there are $(L(G)+o(1))N^k$ edges in the blowup.
As we will mostly talk about the density of a hypergraph rather than the number of edges, it is convenient to define $b(G)=k!L(G)$ to be the \emph{blowup density} of $G$.
Intuitively, it is the largest edge density of the blowups of $G$.
As it turns out, hypergraph Lagrangian is indeed useful for some hypergraph Tur\'an problems, and we will discuss some of those later in \cref{subsec:hypergraph-turan-density} and \cref{sec:known}.

A lesser-known but nonetheless interesting algebraic argument was discovered by Li and Li \cite{LL81}. 
There, they considered the polynomial 
\[f\left((x_v)_{v\in V(G)}\right) = \prod_{uv\not\in E}(x_u-x_v)\]
for any graph $G$.
The key observation is that if $G$ is $K_{r+1}$-free, then $f$ vanishes whenever $r+1$ of the variables $(x_v)_{v\in V(G)}$ are equal to one another.
In light of this, let $I$ be the ideal of polynomials that vanish whenever $r+1$ of the variables are equal.
Then $f\in I$, and Tur\'an's theorem follows from an explicit description of the generators of $I$ that Li and Li worked out.

Our proof looks different from all the proofs mentioned above.
For graphs, our proof can be seen as a double-counting argument that, peculiarly, counts infinitely many objects.
In particular, we will lower bound the number of stars of each size, and show that $K_{r+1}$-freeness actually imposes an upper bound on the numbers.
An interesting feature our proof has is that in order to get the tight bound on the Tur\'an density, it is necessary to take stars of any size into account.
Despite the distinctive look of our proof, our proof is closely related to the standard probabilistic proof of the Caro--Wei theorem.
In fact, if one runs the standard proof on the blowup of the graph, and take the size of the blowup to infinity, then the limit of the argument becomes our argument (we thank Maya Sankar for pointing this out to us).

In spite of the similarity to the proof of the Caro--Wei theorem, our counting argument has the advantage that it can be easily rephrased in terms of entropy.
This will be crucial as it will inform us how we should adapt the proof for hypergraphs.
We will therefore give an introduction to the entropy method in the next subsection.

\subsection{The entropy method}
The concept of entropy in the context of information theory was first formulated by Shannon in his seminal work in 1948 on the noisy-channel coding theorem \cite{Sha48}.
Roughly speaking, the entropy of a random variable measures how much information the random variable carries.
Using entropy, Shannon determined the best efficiency of a code transmitted through a noisy channel that can be corrected with high probability.
This has become the foundation of information theory, and many other definitions of entropy have been made as well.
However, in this paper, we will only use Shannon's definition of entropy.

The adaptation of Shannon entropy in combinatorics and outside the context of information theory came much later in comparison.
Some early examples include Chung, Frankl, Graham and Shearer's work on triangle-intersecting families of graphs \cite{CGFS86} (where Shearer's inequality was introduced), Radhakrishnan's entropic proof of the Bregman's theorem \cite{Rad97}, and Friedgut and Kahn's theorem on the number of copies of a fixed hypergraph in another hypergraph with a given number of edges \cite{FK98}.
There is nonetheless a significant growth in work using the entropy method in the past decade or two.
Two recent exciting, and perhaps unexpected, examples are Gilmer's breakthrough on the union-closed set conjecture \cite{Gil22} and the work of Gowers, Green, Manners and Tao resolving Marton's conjecture (also known as the polynomial Freiman--Ruzsa conjecture over $\FF_2$) \cite{GGMT24}.

In the context of extremal graph theory, the entropy method is particularly useful when dealing with counts of homomorphisms or homomorphism densities.
Here, for any $F,G$ that are graphs or general $k$-graphs, a \emph{homomorphism} from $F$ to $G$ is a function $f:V(F)\to V(G)$ that sends edges of $F$ to edges of $G$. 
In particular, $f$ must be injective on any edge of $F$.
The \emph{homomorphism density} $t(F,G)$ is the probability that a uniformly random chosen function from $V(F)\to V(G)$ is actually a homomorphism.
In this terminology, a corollary of the Kruskal--Katona theorem \cite{Kat68,Kru63} says that $t(K_3, G)\leq t(K_2, G)^{\frac{3}{2}}$, which follows immediately from Shearer's inequality (see also \cite{CY24} for an entropic proof of a slightly stronger result).
In the last decade, the entropy method has been applied to show that various bipartite graphs $F$ are \emph{Sidorenko}, i.e.\ $t(F,G)\geq t(K_2,G)^{e(F)}$.
The prototype of the idea first appeared in the work of Li and Szegedy \cite{LS11}, which was generalized and formulated in terms of entropy by Szegedy \cite{Sze15} and Conlon, Kim, Lee and Lee \cite{CKLL18-1}.
We also refer the readers to \cite{Par14, CL17, CKLL18-2} for related works and \cite{Fitch18, Lee21, GLLV22, BMN24} for related problems.
In fact, in our entropic proofs, we will also derive some Sidorenko-type result using similar arguments.

Given how much the entropy method has been utilized to understand relations between homomorphism densities, it should be surprising that no entropic proof for Tur\'an's theorem was known.
Indeed, an equivalent formulation of the density Tur\'an theorem is that if $t(K_{r+1},G)=0$ then $t(K_2, G)\leq 1-\frac{1}{r}$.
In this paper, we give the first entropic proof of the density Tur\'an theorem.
To do so, we rephrase the density Tur\'an theorem in the following way, and we will later show the equivalence between the two formulations.
Below, and throughout the paper, we use $\HH(X)$ to denote the Shannon entropy of a random variable $X$ (see \cref{sec:prelim} for definitions and basic properties).

\begin{theorem}[Entropic Tur\'an theorem]\label{thm:entropic-turan}
    Let $r$ be a positive integer, and let $G$ be a $K_{r+1}$-free graph.
    Let $X,Y$ be random variables distributed on $V(G)$ so that $\{X,Y\}$ is always an edge in $G$. Assume $X,Y$ are symmetric, i.e.\ the distribution of $(X,Y)$ and the one of $(Y,X)$ are the same.
    Then
    \[\HH(X,Y) \leq 2\HH(X)+\log_2\left(1-\frac{1}{r}\right).\]
\end{theorem}

We make a brief remark that it is easier to see the entropic Tur\'an theorem implies the density Tur\'an theorem by sampling $(X,Y)$ uniformly at random.
To show the equivalence, we use an entropic reinterpretation of blowup density and Langrangian.
Indeed, it turns out that for a given graph $G$, the maximum of the quantity $\HH(X,Y)-2\HH(X)$ for symmetric $V(G)$-valued random variables $X,Y$ with $\{X,Y\}\in E(G)$ is related to the blowup density $b(G)$ of $G$.
More surprisingly, the maximum of $\HH(X,Y)-\HH(X)$ is related to the spectral radius $\rho(G)$ of $G$.
Here, the spectral radius is the largest eigenvalue of the adjacency matrix of $G$.
Those connections will be made precise and proven in \cref{sec:connection}, where we also generalize the connections to hypergraphs.
One benefit is that as an immediate corollary of our entropic Tur\'an theorem, we can generalize spectral Tur\'an theorems established by Wilf \cite{Wil86} and Nikiforov \cite{Nik02,Nik06}.

\begin{theorem}\label{thm:spectral-Turan-tree}
    Let $r\geq 2$ and $T$ be a tree with $\ell\geq 1$ vertices. For any $K_{r+1}$-free graph $G$, we have
    \[\rho(G)^\ell\leq \left(1-\frac{1}{r}\right)\#\{\text{homomorphisms from $T$ to $G$}\}.\]
\end{theorem}
To see that this is indeed a generalization of Wilf's and Nikiforov's results, we can take $T$ to be the path $P_{\ell}$ on $\ell$ vertices.
Wilf's result corresponds to $\ell=1$, whereas Nikiforov's results correspond to $\ell=2$ and general $\ell$. 
\begin{theorem}[\cite{Wil86,Nik02,Nik06}]\label{thm:spectral-Turan}
    Let $r\geq 2$. For any $K_{r+1}$-free graph $G$ with $n$ vertices and $m$ edges, we have
    \[\rho(G)\leq \left(1-\frac{1}{r}\right)n,\]
    \[\rho(G)^2\leq \left(1-\frac{1}{r}\right)\cdot 2m,\]
    and
    \[\rho(G)^\ell\leq \left(1-\frac{1}{r}\right)w_\ell(G),\]
    where $w_\ell(G)$ denotes the number of $\ell$-walks in $G$.
\end{theorem}

\subsection{Hypergraph Tur\'an densities}\label{subsec:hypergraph-turan-density}
Using the idea from our entropic proof of the density Tur\'an theorem, we can determine the Tur\'an densities for some new family of hypergraphs.
Before presenting our results, let us first introduce some definitions and previous work that are relevant.

For any family of $k$-graphs $\cF$, its Tur\'an number $\textup{ex}(n,\cF)$ is defined to be the maximum number of edges in a $k$-graph $G$ that is $F$-free for every $F\in \cF$.
The Tur\'an density is defined analogously by $\pi(\cF) = \lim_{n\to\infty}\textup{ex}(n,\cF)/\binom{n}{k}$.
For any family of $k$-graphs $\cF$ and a $k$-graph $G$, we say that $G$ is \emph{$\cF$-hom-free} if there does not exist any homomorphism $F\to G$ for every $F\in \cF$.
A $F$-hom-free $k$-graph is simply a $k$-graph that is $\{F\}$-hom-free.

To facilitate the discussion, we recall the following standard corollary of supersaturation (see \cite[Section 2]{Kee11} or \cite[Lemma 2.2]{San24} for example).
\begin{theorem}\label{thm:supersaturation}
    For any family $\cF$ of $k$-graphs, its Tur\'an density $\pi(\cF)$ is the supremum of $b(G)$ where $G$ runs through all $\cF$-hom-free $k$-graphs.
\end{theorem}
Notice that a single edge has blowup density $k!/k^k$, showing that $b(G)\geq k!/k^k$ if $G$ is not empty.
This immediately shows that either $\pi(\cF)=0$ or $\pi(\cF)\geq k!/k^k$ for any family of $k$-graphs $\cF$.
We see that among the possible values of Tur\'an density of families of $k$-graphs, there is a ``jump'' going from $0$ to $k!/k^k$.
When $k=2$, this is indeed the behavior of Tur\'an densities: the Erd\H{o}s--Stone theorem shows that all possible values are $0, \frac{1}{2}, \frac{2}{3}, \frac{3}{4},\ldots$, showing that there are only jumps in the case of graphs.
However, for hypergraphs, the set of possible Tur\'an densities has a different behavior.
It was first discovered by Frankl and R\"odl \cite{FR84} that for each $k\geq 3$, there are infinitely many \emph{non-jumps} $\delta$, where for every $\varepsilon>0$ there exists a family $\cF$ of $k$-graphs with $\pi(\cF)\in (\delta,\delta+\varepsilon)$.
On the other hand, Baber and Talbot \cite{BT11} showed that jumps do exist above $k!/k^k$ when $k=3$.
However, our understanding in jumps and non-jumps is still limited, and we do not even know whether $k!/k^k$ is a jump.

A standard argument shows that $k!/k^k$ is a jump if and only if there exists a finite family $\cF$ of $k$-graph with $\pi(\cF)=k!/k^k$ and $b(F)>k!/k^k$ for each $F\in \cF$ (see \cite{FR84}).
The fact that we do not know whether $k!/k^k$ is a jump can thus be seen as a result of not having sufficient understanding in the families $\cF$ with $\pi(\cF)=k!/k^k$.
Indeed, known families with Tur\'an densities equal to $k!/k^k$ are so few that we can list them here.
For general $k$, Mubayi \cite{Mub06} showed that the $k$-uniform extended clique $E^{(k)}_{k+1}$ of size $k+1$ has Tur\'an density $k!/k^k$.
Here, the \emph{extension} of a hypergraph is another hypergraph with higher uniformity obtained by adding different vertices into the edges, and an \emph{extended clique} is an extension of a complete graph. 
In particular, $E^{(k)}_{k+1}$ is obtained by adding $k-2$ extra vertices to each edge of $K_{k+1}$, where no two edges share any extra vertices.
This was later generalized by Mubayi and Pikhurko \cite{MP07}, who showed that the hypergraph $\Delta_{(1,1,\ldots, 1)}$ with edges
\[\left\{v_1,\ldots, v_k\right\}\text{ and }\{w,v_i,u^{(i)}_1,\ldots, u^{(i)}_{k-2}\}\text{ for }i\in [k]\]
also has Tur\'an density $k!/k^k$.
Here, and later whenever the vertex set is not explicitly described, the vertex set consists of vertices that appear in the description of the edges.
Mubayi and Pikhurko's result is indeed an improvement as $E^{(k)}_{k+1}$ is homomorphic to $\Delta_{(1,1,\ldots, 1)}$, showing that $E^{(k)}_{k+1}$-hom-free graphs are also $\Delta_{(1,1,\ldots,1)}$-hom-free and so $\pi(E^{(k)}_{k+1})\leq \pi(\Delta_{(1,1,\ldots,1)})$.

We remark that both Mubayi's \cite{Mub06} and Mubayi and Pikhurko's \cite{MP07} results are stronger---the exact Tur\'an numbers were determined for sufficiently many vertices.
If we only care about the Tur\'an density, then an argument of Sidorenko \cite{Sid89} based on hypergraph Lagrangian can be modified to show that $\pi(\Delta_{(1,\ldots,1)})=k!/k^k$ as well---this is an observation by Keevash \cite[Theorem 3.1]{Kee11}.

For smaller $k$'s, slightly more is known.
When $k=3$, Bollob\'as \cite{Bol74} showed that $\pi(\{K_4^{-},F_5\}) = \frac{2}{9}$ where $K_4^{-} = \{123,124,134\}$ and $F_5=\{123,124,345\}$.
This was improved by Frankl and F\"uredi \cite{FF83}, who showed that $\pi(F_5)$ is already equal to $\frac{2}{9}$.
Using flag algebra, Baber and Talbot \cite{BT12} improved this further by showing that $\pi(\{123,124,345,156\}) = \frac{2}{9}$.
Finally, when $k=4$, Pikhurko \cite{Pik08} showed that $\pi(\{1234, 1235, 4567\}) = \frac{3}{32}$.

As shown above, not a lot is known about families $\cF$ of $k$-graphs with $\pi(\cF)=k!/k^k$.
As an application of our entropic proof of the density Tur\'an theorem, we will generalize our argument to show $\pi(\cF)=k!/k^k$ for a new family $\cF$ of $k$-graphs.
Our method has a benefit that we may first come up with an argument and then see what family of $k$-graphs need to be forbidden in order for the argument to work.
We believe that this advantage can help discovering more families $\cF$ with minimum positive Tur\'an densities.

\begin{figure}[h]\centering\label{fig:Tent}
\begin{tikzpicture}[scale=0.8]

    \coordinate (A) at (0,0);
    \coordinate (B) at (1,0);
    \coordinate (C) at (2,0);
    \coordinate (D) at (5,0);
    \coordinate (E) at (6,0);
    \coordinate (F) at (4.5,1.732/2);
    \coordinate (G) at (2.75,1.732*3/4);
    \coordinate (H) at (4,1.732);
    \coordinate (I) at (3.5,1.732*3/2);

    \draw [fill] (A) circle (1.6pt);
    \draw [fill] (B) circle (1.6pt);
    \draw [fill] (C) circle (1.6pt);
    \draw [fill] (D) circle (1.6pt);
    \draw [fill] (E) circle (1.6pt);
    \draw [fill] (F) circle (1.6pt);
    \draw [fill] (G) circle (1.6pt);
   \draw [fill] (H) circle (1.6pt);
    \draw [fill] (I) circle (1.6pt);
% Define a path with rounded corners
\draw[rounded corners=8pt,black,line width=2pt] 
    (0-0.2,0.5)--(-0.5-0.2,0)--(0-0.2,-0.5)--(6+0.2,-0.5)--(6.5+0.2,0)--(6+0.2,0.5)--cycle;
    
\draw[rounded corners=6pt,black,line width=2pt] 
(0-0.2,0.3)--(-0.3-0.2,0)--(0-0.2,-0.3)--(2+0.3/1.732,-0.3)--(3.5+0.1+0.15*1.732,1.732*3/2+0.1*1.732-0.15)--(3.5+0.1+0.15,1.732*3/2+0.1*1.732+0.15*1.732)--(3.5+0.1-0.15*1.732,1.732*3/2+0.1*1.732+0.15)--(2-0.3/1.732,0.3)--cycle;

\draw[rounded corners=6pt,black,line width=2pt] 
(6+0.2,0.3)--(6+0.3+0.2,0)--(6+0.2,-0.3)--(5-0.3/1.732,-0.3)--(3.5-0.1-0.15*1.732,1.732*3/2+0.1*1.732-0.15)--(3.5-0.1-0.15,1.732*3/2+0.1*1.732+0.15*1.732)--(3.5-0.1+0.15*1.732,1.732*3/2+0.1*1.732+0.15)--(5+0.3/1.732,0.3)--cycle;

\node at (3.5,0) {Base};
\node at (4.7,1.732*3/2) {Apex};
\end{tikzpicture}
\caption{$(3,2)$-tent}
\end{figure}
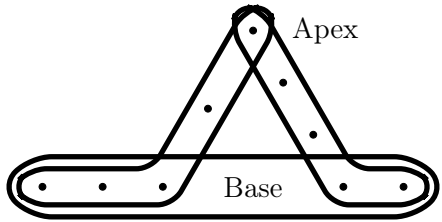

To state our result, for any partition $\lambda$ of $k$, let $\lambda = (\lambda_1,\ldots, \lambda_{\ell})$ where $\ell = \ell(\lambda)$ is the length of $\lambda$, and $\lambda_1\geq \cdots\geq \lambda_{\ell}$.
We also denote $\sum_{i=1}^{\ell}\lambda_i$ by $\abs{\lambda}$ (which is equal to $k$ by definition).
For any $\lambda$ with $\ell(\lambda)\geq 2$, we define the \emph{$\lambda$-tent}, denoted by $\Delta_{\lambda}$, to be the following $k$-graph.
The $\lambda$-tent comes with an edge $e$ that is the \emph{base} and a vertex $v$ that is the \emph{apex}.
Setting $\ell=\ell(\lambda)$ to be the length of $\lambda$, for each $i\in[\ell]$ we also have an edge $e_i$ containing $v$ such that $\abs{e_i\cap e}=\lambda_i$.
Moreover, we require that $e_i\cap e_j = \{v\}$ for any $i\neq j\in [\ell]$.
It is clear that this determines $\Delta_{\lambda}$ uniquely up to isomorphism---in fact, we must have $e\cap e_1,\ldots, e\cap e_{\ell}$ partition $e$.
It is easy to check that this definition matches the definition of $\Delta_{(1,1,\ldots,1)}$ above, that $F_5 = \Delta_{(2,1)}$ (with base $123$ and $4$ being the apex), and that Pikhurko's result can be rephrased as $\pi(\Delta_{(3,1)})=\frac{3}{32}$.
Our result can now be stated as follows.

\begin{theorem}\label{thm:main-tent}
    Let $k\geq 2$ be a positive integer, and let $\cF_k$ be the family of $\lambda$-tents with $\abs{\lambda}=k$ and $\ell(\lambda)=2$.
    Then $\pi(\cF_k) = k!/k^k$.
\end{theorem}

Note that this is a stronger statement than Mubayi's and Mubayi and Pikhurko's results.
In fact, $\Delta_{(1,1,\ldots, 1)}$ admits a homomorphism to $\Delta_{\lambda}$ for every $\abs{\lambda}=k$ and $\ell(\lambda)=2$, which shows that $\pi(\Delta_{(1,1,\ldots,1)})\leq \pi(\cF_k)$. 
Using the same argument, we can transform \cref{thm:main-tent} into a Tur\'an result of a single $k$-graph.

\begin{theorem}\label{thm:one-tent}
    Let $k\geq 2$ be a positive integer, and let $\lambda$ be a partition of $k$ such that $\lambda_1\leq \lceil k/2\rceil$ and $\lambda_i=1$ for all $1<i\leq \ell(\lambda)$.
    Then $\pi(\Delta_\lambda) = k!/k^k$.
\end{theorem}

Although when $k=3$ and $4$, \cref{thm:one-tent} is subsumed by the known results mentioned above, this gives a new Tur\'an result for larger $k$'s. 
To show that this should be a nontrivial result for larger $k$'s, we prove the following result in the opposite direction.

\begin{theorem}\label{thm:not-Turan}
    There exists a constant $\alpha<1$ so that for all sufficiently large $k\in\NN$ and any partition $\lambda$ of $k$ with $\ell(\lambda)\geq 2$, if $\lambda_1 >\alpha k$ then $\pi(\Delta_{\lambda})>k!/k^k$.
\end{theorem}

\cref{thm:one-tent} shows that the constant in \cref{thm:not-Turan} cannot be smaller than $1/2$, and it seems like an interesting question to determine the best possible value of $\alpha$.
It might help us understand the $k$-graphs $F$ with $\pi(F)=k!/k^k$ as well.
We leave this as a future direction for interested readers.

Beyond showing $\pi(\cF) = k!/k^k$ for various families $\cF$ of $k$-graphs, our method also applies to some other scenarios where the extremizers are blowups of complete hypergraphs.
Unfortunately, we have not been able to find an argument that proves a new and clean statement in those settings.
We nonetheless include the arguments later in \cref{sec:known} in the hope that they will be enlightening for readers interested in adapting our arguments.
The relevant background will also be introduced there.

In those proofs, we make no effort in deciding the structures of extremizers, but we suspect that it is possible to extract some information from our proofs.
We briefly discuss this in \cref{subsec:stability}.
\subsection{Structure of the paper}
We will first present our new proof of the density Tur\'an theorem in \cref{sec:counting}.
We will then introduce the necessary entropic tools in \cref{sec:prelim}, which will set us up for \cref{sec:entropy}, where we rephrase our proof in terms of entropy.
In \cref{sec:connection}, we will show how our entropic formulation captures quantities such as hypergraph Lagrangian and spectral radius.
We will use the connection to prove the spectral Tur\'an theorems and the equivalence between the entropic Tur\'an theorem and the density Tur\'an theorem.
In \cref{sec:partial-hypergraph}, we set up some notations and propositions that will be useful in the later sections.
In \cref{sec:main-proof}, we will apply the entropic argument in \cref{sec:entropy} to show \cref{thm:main-tent} in two different ways, and we will also prove \cref{thm:one-tent,thm:not-Turan}.
Some further generalization of our arguments is included in \cref{sec:known}, where we also introduce some related known results.
Finally, we will end with some concluding remarks in \cref{sec:conclusion}.

\section{A new proof of the density Tur\'an theorem}\label{sec:counting}
In this section, we give a new proof to the density Tur\'an theorem. The key idea is to lower bound the density of stars of each size in terms of edge density by their Sidorenko property. If the densities are large, then we shall find a large clique. The main difference of this proof from all the previous ones is that we consider stars of all sizes at once.

\begin{proof}[Proof of the density Tur\'an theorem]
    For any two graphs $H,G$, let $t(H,G)$ be the homomorphism density of $H$ in $G$. That is, $t(H,G)$ is the probability that a function $f:V(H)\rightarrow V(G)$ chosen uniformly at random is a homomorphism from $H$ to $G$. We will need the following lemma about lower bounding the homomorphism density of stars in terms of edge density, which is a special case of Sidorenko's conjecture. We include the proof here since the proof is short.
    \begin{lemma}
        For $i\geq 0$, let $S_i=K_{1,i}$ be the star with $i+1$ vertices. Then
        \[t(S_i,G)\geq t(K_2,G)^i\]
        holds for any graph $G$.
    \end{lemma}
    \begin{proof}
        Assume $n=\abs{V(G)}$ and $m=\abs{E(G)}$. Note that $S_i$ has $i+1$ vertices, and hence 
        \[t(S_i,G)=\frac{\sum_{v\in V(G)}\deg(v)^i}{n^{i+1}}\geq \frac{1}{n^i}\left(\frac{\sum_{v\in V(G)}\deg(v)}{n}\right)^i=\frac{(2m)^i}{n^{2i}}=t(K_2,G)^i,\]
        where the inequality follows from the convexity of $x^i$.
    \end{proof}
    Now we assume the graph $G$ is $K_{r+1}$-free. We sample a sequence of i.i.d.\@ random vertices $v_0,v_1,\dots$ from $V(G)$ uniformly at random. For $i\geq 0$, let $A_i$ be the event that the induced graph on vertices $v_0,\dots,v_{i-1},v_i$ contains $S_i$ as a subgraph centered at $v_i$. In particular, $A_0$ is the true event. Note that there can only be at most $r$ events happening at the same time. Otherwise, assume $A_{i_0},A_{i_1},\dots,A_{i_r}$ are all true for some $0=i_0<i_1<\dots<i_r$. Then $v_{i_0},\dots,v_{i_r}$ form an $(r+1)$-clique in $G$. Therefore, by double counting, we may conclude that
    \[\PP(A_0)+\PP(A_1)+\dots\leq r.\]
    
    On the other hand, we know that $\PP(A_i)=t(S_i,G)\geq t(K_2,G)^i$ for all $i$. Thus, we have
    \[\frac{1}{1-t(K_2,G)}\leq \PP(A_0)+\PP(A_1)+\dots\leq r.\]
    After rearranging, we get
    \[\frac{2m}{n^2}=t(K_2,G)\leq 1-\frac{1}{r},\]
    and we are done.
\end{proof}

\section{Shannon entropy}\label{sec:prelim}
In this section, we introduce the definition of Shannon entropy and some of the properties we will use from the literature. We refer the readers to \cite[Section 14.6]{AS00} for a more detailed introduction.
We will also prove a lemma which upper bounds the entropies of random variables by the entropy of their mixture. This lemma will be one of the key ingredients of many of the proofs in the rest of this paper.

\subsection{Preliminaries}

For any discrete random variable $X$, we write $p_X(x)\eqdef\PP(X=x)$. Also, we denote by $\supp(X)$ the support of $X$, i.e.\ the set of all $x$ such that $p_X(x)>0$.
Throughout this paper, the random variables we will consider are always discrete with finite support, i.e.\ $\abs{\supp(X)}<\infty$. For any such random variable, we may define its Shannon entropy.

\begin{definition}
    For any random variable $X$, we define its Shannon entropy
    \[\HH(X)\eqdef\sum_{x\in\supp(X)}-p_X(x)\log_2p_X(x).\]
    For any sequence of random variables $X_1,\dots,X_n$, we use $\HH(X_1,\dots,X_n)$ to denote the entropy of the random tuple $(X_1,\dots,X_n)$.
\end{definition}

We also define the conditional entropy of $X$ given $Y$.
\begin{definition}
    For any two random variables $X,Y$, the conditional entropy of $X$ given $Y$ is given by
    \[\HH(X\mid Y)\eqdef\HH(X,Y)-\HH(Y).\]
    Equivalently, we have
    \begin{align*}
        \HH(X\mid Y)=&\quad\smashoperator{\sum_{y\in\supp(Y)}}p_Y(y)\HH(X\mid Y=y)\\
        =&\quad\smashoperator{\sum_{(x,y)\in \supp(X,Y)}}-p_{X,Y}(x,y)\log_2\left(\frac{p_{X,Y}(x,y)}{p_Y(y)}\right).
    \end{align*}
\end{definition}

Using the definition of conditional entropy, we have the following chain rule.
\begin{proposition}[Chain rule]
    For any random variables $X_1,\dots,X_n$, we have
    \[\HH(X_1,\dots,X_n)=\HH(X_1)+\HH(X_2\mid X_1)+\dots+\HH(X_n\mid X_1,\dots,X_{n-1}).\]
\end{proposition}

The following proposition says that on a fixed support, the entropy is maximized by the uniform distribution on that support.
\begin{proposition}[Uniform bound]
    For any random variable $X$, we have
    \[\HH(X)\leq \log_2\abs{\supp(X)},\]
    where the equality holds if and only if $X$ is uniform.
\end{proposition}

We will also need the following two propositions about entropy.
\begin{proposition}[Subadditivity]
For any three random variables $X,Y,Z$, we have
\[\HH(X,Y)\leq \HH(X)+\HH(Y),\]
\[\HH(X,Y\mid Z)\leq \HH(X\mid Z)+\HH(Y\mid Z).\]
\end{proposition}
\begin{proposition}[Dropping condition]
For any three random variables $X,Y,Z$, we have
\[\HH(X\mid Y)\leq \HH(X),\]
\[\HH(X\mid Y,Z)\leq \HH(X\mid Z).\]
\end{proposition}

\subsection{Mixture and the mixture bound}
In this subsection, the concern is what is called the \emph{mixture} of random variables.
\begin{definition}
    For random variables $X_1,\dots,X_n$ and weights $w_1,\dots,w_n\geq 0$ with $\sum_{i=1}^n w_i=1$, we say that $Z$ is the \emph{mixture of $X_1,\dots,X_n$ with weight $w_1,\dots, w_n$} if $Z$ is obtained from the following procedure. We first pick an independent random index $\bi$ with probability $\PP(\bi =i)=w_i$. Then we set $Z=X_{\bi}$.
\end{definition}
In our applications, we will consider mixtures of random variables whose supports do not overlap too much.
\begin{definition}
    Let $a$ be a positive integer.
    We say that the random variables $X_1,\dots,X_n$ have \emph{$(a+1)$-wise disjoint supports} if for any element $x\in \cup_{i=1}^n \supp (X_i)$, there are at most $a$ many indices $i$ such that $x\in\supp (X_i)$.
\end{definition}

With the definitions above, we may state our lemma about an upper bound on the entropies of random variables with $(a+1)$-wise disjoint supports, in terms of the entropy of their mixture.
\begin{lemma}[Mixture bound]\label{lemma:mix}
    Let $X_1,\dots,X_n$ be random variables with $(a+1)$-wise disjoint supports. Then there exists a mixture of $X_1,\dots,X_n$, say $Z$, such that
    \[\sum_{i=1}^n 2^{\HH(X_i)}\leq a2^{\HH(Z)}.\]
\end{lemma}
Before proving the lemma, we use the following example to illustrate that \cref{lemma:mix} resembles a double counting on $(a+1)$-wise disjoint sets.
\begin{example}
    Let $a$ be an integer and let $S_1,\dots,S_n$ be some sets that are $(a+1)$-wise disjoint. Assume $X_i$ is a random element chosen from $S_i$ uniform at random for each $i\in [n]$, and let $Z$ be the mixture of $X_1,\ldots, X_n$ provided by \cref{lemma:mix}. We have $2^{\HH(X_i)}=\abs{S_i}$, and by uniform bound we have $2^{\HH(Z)}\leq \abs{\cup_{i=1}^nS_i}$.
    Hence, \cref{lemma:mix} implies that
    \[\sum_{i=1}^n\abs{S_i}\leq a2^{\HH(Z)}\leq a\abs{\bigcup_{i=1}^n S_i},\]
    which gives the same bound as the double counting argument on pairs $(x,i)$ with $x\in S_i$.
\end{example}
Thus, the mixture bound can be viewed as an entropic version of this double counting.

\begin{proof}[Proof of \cref{lemma:mix}]
    Let $s_i=2^{\HH(X_i)}$ and we define
    \[w_i=\frac{s_i}{\sum_{j=1}^n s_j}.\]
    Let $\bi$ be an independent random index with probability $\PP(\bi =i)=w_i$ and let $Z=X_{\bi}$ be the mixture. By the chain rule, we have $\HH(Z,\bi)=\HH(\bi)+\HH(Z\mid \bi)=\HH(Z)+\HH(\bi\mid Z)$. Therefore,
    \[\HH(Z)=\HH(\bi)+\HH(Z\mid \bi)-\HH(\bi\mid Z).\]
    By the definition of entropy and conditional entropy, we have
    \[\HH(\bi)=\sum_{i=1}^n -w_i\log_2 w_i=\sum_{i=1}^n\frac{-s_i}{\sum_{j=1}^n s_j}\log_2\bigl(\frac{s_i}{\sum_{j=1}^n s_j}\bigr)\]
    and
    \[\HH(Z\mid \bi)=\sum_{i=1}^n w_i\HH(X_i)=\sum_{i=1}^n\frac{s_i\log_2 s_i}{\sum_{j=1}^n s_j}.\]
    We may upper bound $\HH(\bi\mid Z)$ by uniform bound. For any $x\in\cup_{i=1}^n\supp(X_i)$, when conditioning on $Z=x$, there are at most $a$ possible indices as an outcome of $\bi$. Thus, we have
    \[\HH(\bi\mid Z)\leq \log_2 a.\]
    Combining all above, we get
    \begin{align*}
        \HH(Z)\geq &\sum_{i=1}^n\left(\frac{-s_i}{\sum_{j=1}^n s_j}\log_2\bigl(\frac{s_i}{\sum_{j=1}^n s_j}\bigr)+\frac{s_i\log_2 s_i}{\sum_{j=1}^n s_j}\right)-\log_2 a\\
        =&\log_2\left(\sum_{j=1}^n s_j\right)-\log_2 a,
    \end{align*}
    and we are done after rearranging.
\end{proof}

\section{Reformulation using the entropy method}\label{sec:entropy}
In this subsection, we reformulate the proof in \cref{sec:counting} using entropy to prove \cref{thm:entropic-turan}. As expected, we shall sample the stars in the same way as in \cite{Sze15,CKLL18-1}, and we will use \cref{lemma:mix} to replace the double counting argument.
\begin{proof}[Proof of \cref{thm:entropic-turan}]
    Recall that we have a $K_{r+1}$-free graph $G$ and symmetric random variables $X,Y$ distributed on $V(G)$ with $\{X,Y\}\in E(G)$ always holding.
    We first fix an integer $N\in\NN$, and we will take $N$ goes to infinity later. 
    \begin{claim}\label{claim:StarSido}
        For each $i=0,1,\dots, N$, there exists a random tuple  $T_i=(v_0^{(i)},\dots,v_N^{(i)})\in V(G)^{N+1}$ such that 
        \begin{enumerate}
            \item there is always an edge between $v_j^{(i)},v_i^{(i)}$ for all $j=0,\dots,i-1$,
            \item the marginal distributions of $v_j^{(i)}$ and $X$ are the same for all $j=0,1\dots,N$, and
            \item $\HH(T_i)=i\HH(Y\mid X)+(N+1-i)\HH(X)$.
        \end{enumerate}
    \end{claim}
    \begin{proof}
        For $i=0$, it is easy to check that $N+1$ i.i.d. random vertices $v_0^{(0)},\dots,v_N^{(0)}$ with the law of $X$ satisfy the condition.

        For $i\geq 1$, we first sample an edge $(v_0^{(i)},v_i^{(i)})$ using the law of $(X,Y)$. Next, we condition on $v_i^{(i)}$ and resample $v_0^{(i)}$ $(i-1)$ times conditionally independently to get $v_1^{(i)},\dots,v_{i-1}^{(i)}$. Finally, we sample $v_{i+1}^{(i)},\dots, v_N^{(i)}$ independently using the law of $X$. 

        Note that the first two conditions are true from the way we sample the random variables. It remains to compute $\HH(T_i)$. Note that $\HH(T_i)=\HH(v_0^{(i)},\dots,v_i^{(i)})+(N-i)\HH(X)$ since we sampled $v_{i+1}^{(i)},\dots, v_N^{(i)}$ independently. By the chain rule, we have 
        \begin{align*}            \HH(v_0^{(i)},\dots,v_i^{(i)})=&\HH(v_0^{(i)},\dots,v_{i-1}^{(i)}\mid v_i^{(i)})+\HH(v_i^{(i)})\\
        =&i\HH(v_0^{(i)}\mid v_i^{(i)})+\HH(v_i^{(i)})\\
        =&i\HH(Y\mid X)+\HH(X).
        \end{align*}
        Therefore, $\HH(T_i)=i\HH(Y\mid X)+(N+1-i)\HH(X)$.
    \end{proof}
    Now, we may apply \cref{lemma:mix} to the random tuples $T_0,\dots,T_N$ in \cref{claim:StarSido}. Since $G$ is $K_{r+1}$-free, similar to the proof in \Cref{sec:counting}, any tuple of $N+1$ vertices is in at most $r$ supports $\supp(T_i)$. Therefore, the supports of $T_0,\dots,T_N$ are $(r+1)$-wise disjoint. Thus, there is a mixture $T=(v_0,\dots,v_N)$ of $T_0,\dots,T_N$ such that
    \[\sum_{i=0}^N 2^{\HH(T_i)}\leq r2^{\HH(T)}.\]
    Note that the marginal distribution of $v_i$ is also the same as the marginal distribution of $X$, so we may upper bound $\HH(T)$ by $(N+1)\HH(X)$ by subadditivity. By using $\HH(T_i)=i\HH(Y\mid X)+(N+1-i)\HH(X)$, we get
    \[\sum_{i=0}^N x^i\leq r,\]
    where $x\eqdef 2^{\HH(Y\mid X)-\HH(X)}$. By taking $N$ to infinity, we conclude that $1/(1-x)\leq r$. Therefore,
    \[\HH(Y\mid X)-\HH(X)=\log_2x\leq \log_2\left(1-\frac{1}{r}\right).\qedhere\]
\end{proof}

    Let $\abs{V(G)}=n$ and $\abs{E(G)}=m$. 
    If we pick $(X,Y)$ uniformly at random from all the oriented edges, \cref{thm:entropic-turan} and the uniform bound give
    \[\log_2(2m)=\HH(X,Y)\leq 2\HH(X)+\log_2\left(1-\frac{1}{r}\right)\leq 2\log_2 n+\log_2\left(1-\frac{1}{r}\right).\]
    That is, $m\leq \left(1-\frac{1}{r}\right)\frac{n^2}{2}$, which recovers the density Tur\'an theorem.
    In the next section, we will see that \cref{thm:entropic-turan} is in fact equivalent to the density Tur\'an theorem by relating entropy to blowup densities.

\section{Connecting entropy to Lagrangian and spectral radius}\label{sec:connection}

In this section, we will show that \cref{thm:entropic-turan} is equivalent to the density Tur\'an theorem.
We will actually generalize this equivalence in many ways: we will show it for hypergraphs, and we will also go much beyond Lagrangian and blowup densities.
This will be useful later to draw connection to the spectral radius of graphs.

We first observe that in \cref{thm:entropic-turan}, the quantity that we care about is actually the maximum of $\HH(X,Y)-2\HH(X)$ when $(X,Y)$ ranges over all possible symmetric distributions on the oriented edges of $G$.
This quantity turns out to be related to the blowup density $b(G)$.
To extend this to hypergraphs, we make the following definitions.

\begin{definition}[Random edge with uniform ordering]
    Let $G$ be a $k$-graph, we say that a tuple of random vertices $(X_1,\dots,X_k)\in V(G)^k$ is a \emph{random edge with uniform ordering on $G$} if $(X_1,\dots,X_k)$ is symmetric and $\{X_1,\dots,X_k\}$ is always an edge of $G$. Here, $(X_1,\dots,X_k)$ being symmetric means the distribution of $(X_{\sigma(1)},\dots,X_{\sigma(k)})$ is always the same for any permutation $\sigma$ of $[n]$.
\end{definition}

\begin{definition}[Entropic density]
    
    For any $k$-graph $G$, define its \emph{entropic density} $b_{\textup{entropy}}(G)$ to be the largest possible value of $2^{\HH(X_1,\ldots, X_k)-k\HH(X_1)}$ for any random edge with uniform ordering $(X_1,\ldots, X_k)$.
\end{definition}

Note that $b_{\textup{entropy}}(G)$ exists as the space of random edge with uniform ordering is compact.
We will show that $b_{\textup{entropy}}(G)$ is equal to $b(G)$, which immediately shows that \cref{thm:entropic-turan} is equivalent to the density Tur\'an theorem.
We will actually show a stronger statement.
To that end, we make the following notations.
For any $k$-graph $G$, let $\vec{E}(G)$ be the set of oriented edges, i.e.\ tuples $(v_1,\ldots, v_k)\in V(G)^k$ with $\{v_1,\ldots,v_k\}\in E(G)$.
For each $p>0$, let $b_p(G)$ be the maximum of $\prod_{(v_1,\ldots, v_k)\in\vec{E}(G)}x_{v_1}\cdots x_{v_k}$ for $(x_v)_{v\in V(G)}$ subject to $\norm{x_v}_{\ell^p}=1$ (the same definition was made by Keevash, Lenz and Mubayi \cite{KLM14} where they called the quantity the \emph{$p$-spectral radius}).
Also let $b_{p,\textup{entropy}}(G)$ be the largest possible value of $2^{\HH(X_1,\ldots, X_k)-\frac{k}{p}\HH(X_1)}$ for any random edge with uniform ordering $(X_1,\ldots, X_k)$.
Note that $b_p(G)$ and $b_{p,\textup{entropy}}(G)$ both exist by compactness.

\begin{example}\label{ex:entropic-lagrangian}
    When $p=1$, we clearly have $b_p(G)=b(G)$ and $b_{p,\textup{entropy}}(G) = b_{\textup{entropy}}(G)$.
    When $G$ is a graph and $p=2$, it is not hard to see that $b_p(G)$ is the maximum 
    \[\max\vec{x}^{\intercal}A_G\vec{x}\textup{ subject to }\norm{(x_v)_{v\in V(G)}}_{\ell^2}=1\]
    where $A_G$ is the adjacency matrix of $G$.
    It is a standard fact that this is exactly the spectral radius of $G$.
    In this case, $b_{2,\textup{entropy}}(G)$ is the largest possible value of $2^{\HH(X,Y)-\HH(X)} = 2^{\HH(Y\mid X)}$ for any random edge with uniform ordering $(X,Y)$. 

    For general $k$, if $p=k$, then $b_p(G)$ corresponds to the spectral radius of the adjacency \linebreak$k$-tensor of $G$, which was proven in \cite{Qi13}.
    The quantity $b_{k,\textup{entropy}}(G)$ is the largest possible value of $2^{\HH(X_1,\ldots, X_k)-\HH(X_1)} = 2^{\HH(X_2,\ldots, X_k\mid X_1)}$. Once we prove $b_k(G) = b_{k,\textup{entropy}}(G)$, this would provide a nice alternative interpretation of the spectral radius for hypergraphs.
\end{example}

Now we will show that $b_p(G)$ and $b_{p,\textup{entropy}}(G)$ are equal to each other.
The proof uses Lagrange multiplier in a crucial way.

\begin{proposition}\label{prop:entropic-lagrangian}
    For any $k$-graph $G$ and any $p>0$, $b_{p,\textup{entropy}}(G) = b_p(G)$.
    
\end{proposition}

\begin{proof}
    For any $v\in V(G)$, let $\vec{L}_v(G)$ be the oriented link of $v$, i.e.\ the set of $(v_2,\ldots, v_k)$ such that $(v,v_2,\ldots, v_k)\in \vec{E}(G)$.
    
    We start with the following claim that helps us simplify $\HH(X_1,\dots,X_k)-\frac{k}{p}\HH(X_1)$ when $(X_1,\dots,X_k)$ is in a certain form.
    \begin{claim}\label{claim:entropy-lag}
        For any tuple $(x_v)_{v\in V(G)}\in \RR_{\geq 0}^{V(G)}$, we consider a random edge with uniform ordering $(X_1,\dots,X_k)$ on $G$ given by
    \[\PP((X_1,\dots,X_k)=(v_1,\dots,v_k))=\frac{1}{\beta}\prod_{i=1}^kx_{v_i},
    \text{ where }\beta\eqdef\sum_{(v_1,\dots,v_k)\in \vec{E}(G)}\prod_{i=1}^kx_{v_i}.\]
    We also define
    \[y_v\eqdef \PP(X_1=v)=\frac{x_v}{\beta}\sum_{(v_2,\dots,v_k)\in \vec{L}_v(G)}\prod_{i=2}^kx_{v_i}.\]
    Then we have
    \[\HH(X_1,\dots,X_k)-\frac{k}{p}\HH(X_1)=\log_2\beta-\frac{k}{p}\sum_{v\in V(G)}y_v\log_2 \left(\frac{x_v^p}{y_v}\right).\]
    \end{claim}
    \begin{proof}
         First, we have
    \begin{align*}
        \HH(X_1,\dots,X_k)=&\sum_{(v_1,\dots,v_k)\in \vec{E}(G)}-\frac{1}{\beta}\prod_{i=1}^kx_{v_i}\log_2\left(\frac{1}{\beta}\prod_{i=1}^kx_{v_i}\right)\\
        =&\sum_{(v_1,\dots,v_k)\in \vec{E}(G)}\frac{1}{\beta}\prod_{i=1}^kx_{v_i}\left(\log_2\beta-\sum_{i=1}^k\log_2x_{v_i}\right)\\
        =&\log_2\beta-k\sum_{v\in V(G)}y_v\log_2 x_v
    \end{align*}
    Combining this with $\HH(X_1)=\sum_{v\in V(G)}-y_v\log_2y_v$, we get
    \begin{align*}
        \HH(X_1,\dots,X_k)-\frac{k}{p}\HH(X_1)=&\log_2\beta-\frac{k}{p}\sum_{v\in V(G)}\left(py_v\log_2 x_v-y_v\log_2y_v\right)\\
        =&\log_2\beta-\frac{k}{p}\sum_{v\in V(G)}y_v\log_2 \left(\frac{x_v^p}{y_v}\right).\qedhere
    \end{align*}
    \end{proof}
    Now, we may prove the proposition. We first show that $b_{p,\textup{entropy}}(G)\geq b_p(G)$. 
    
    Let $(x_v)_{v\in V(G)}\in \RR_{\geq 0}^{V(G)}$ be the tuple that achieves the maximum in the definition of $b_p(G)$. Define $(X_1,\dots,X_k)$, $\beta$, and $(y_v)_{v\in V(G)}$ in the same way as in \cref{claim:entropy-lag}. Note that $\beta=b_p(G)$ and $\sum_{v\in V(G)}x_v^p=1$. From \cref{claim:entropy-lag}, we have
    \begin{align*}
        \HH(X_1,\dots,X_k)-\frac{k}{p}\HH(X_1)=&\log_2\beta-\frac{k}{p}\sum_{v\in V(G)}y_v\log_2 \left(\frac{x_v^p}{y_v}\right)\\
        \geq& \log_2\beta-\frac{k}{p}\log_2 \left(\sum_{v\in V(G)}x_v^p\right)=\log_2\beta,
    \end{align*}
    where the inequality follows from the Jensen's inequality and the concavity of $\log_2x$.
    Therefore $b_{p,\textup{entropy}}(G)\geq b_p(G)$.

    For the opposite direction, let $(X_1,\ldots, X_k)$ be a random edge with uniform ordering achieving the maximum of $b_{p,\textup{entropy}}(G)$.
    For any unoriented edge $e\in E(G)$, let $q_e$ be the probability $\PP(\{X_1,\ldots, X_k\}=e)$.
    Also let $x_v = \left(\frac{1}{k}\sum_{e\ni v} q_e\right)^{1/p}$.
    Then
    \[\HH(X_1,\ldots, X_k) = \HH(X_1,\ldots, X_k\mid \{X_1,\ldots, X_k\})+\HH(\{X_1,\ldots, X_k\}) = \log_2 k! -\sum_{e\in E(G)}q_e\log_2 q_e\]
    and
    \[\HH(X_1) = \sum_{v\in V}-x_v^p\log_2 x_v^p.\]
    Therefore, $(q_e)_{e\in E(G)}$ is a maximizer of 
    \[-\sum_{e\in E(G)}q_e\log_2 q_e+\frac{k}{p}\sum_{v\in V(G)}x_v^p\log_2 x_v^p\]
    subject to $q_e \geq 0$ for all $e\in E(G)$ and $\sum_{e\in E(G)}q_e=1$.
    Note that $\partial x_v^p/\partial q_e$ is nonzero only if $v\in e$, and if that is the case we have $\partial x_v^p/\partial q_e=1/k$.
    By Lagrange multiplier, we know that 
    \[-\log_2 q_e-\frac{1}{\log 2} + \frac{1}{p}\sum_{v\in e}\left(\frac{1}{\log 2}+\log_2 x_v^p\right)\]
    is constant for all $e\in E(G)$ with $q_e>0$.
    Therefore
    \[\alpha \eqdef \frac{q_e}{\prod_{v\in e}x_v}\]
    is the same for all $e\in E(G)$ with $q_e>0$.
    Notice that $\PP(X_1 = v) = x_v^p$ for any $v\in V(G)$, and for any $(v_1,\ldots, v_k)\in \vec{E}(G)$, we have 
    \[\PP((X_1,\ldots, X_k)=(v_1,\ldots, v_k)) = \frac{q_{\{v_1,\dots,v_k\}}}{k!}=\frac{\alpha}{k!} \prod_{i=1}^{k}x_{v_i}.\]
    Therefore, using \Cref{claim:entropy-lag} with $\beta=k!/\alpha$, we see that
    \begin{align*}
        \HH(X_1,\dots,X_k)-\frac{k}{p}\HH(X_1)=&\log_2\beta-\frac{k}{p}\sum_{v\in V(G)}y_v\log_2 \left(\frac{x_v^p}{y_v}\right),
    \end{align*}
    where, in this case, $y_v=x_v^p$. Thus, $\HH(X_1,\dots,X_k)-\frac{k}{p}\HH(X_1)=\log_2 \beta$.
    Note that $\sum_{v\in V(G)}x_v^p=1$. 
    Therefore by the fact that 
    \[\beta = \sum_{(v_1,\ldots, v_k)\in \vec{E}(G)}\prod_{i=1}^{k}x_{v_i},\]
    we have $b_{p,\textup{entropy}}(G)\leq b_p(G)$.
\end{proof}

\begin{corollary}\label{cor:equiv-to-entropy}
    For any family $\cF$ of $k$-graphs, $\pi(\cF)$ is the supremum of $2^{\HH(X_1,\ldots, X_k)-k\HH(X_1)}$ for any random edge with uniform ordering $(X_1,\ldots, X_k)$ on any $\cF$-hom-free $k$-graph $G$.
\end{corollary}
\begin{proof}
    Since $\pi(\cF)$ is the supremum of $b(G)$ for all $\cF$-hom-free $k$-graphs $G$ by \cref{thm:supersaturation}, we know that $\pi(\cF)$ is the supremum of $b_{\textup{entropy}}(G)$ for all $\cF$-hom-free $k$-graphs $G$ as well.
    The statement follows from the definition of entropic density $b_{\textup{entropy}}(G)$.
\end{proof}

\begin{corollary}\label{cor:entropy-equiv-to-Turan}
    The entropic Tur\'an theorem (\cref{thm:entropic-turan}) is equivalent to the density Tur\'an theorem.
\end{corollary}
\begin{proof}
    By \cref{cor:equiv-to-entropy}, it suffices to show that if $G$ is $K_{r+1}$-free, then $G$ is $K_{r+1}$-hom-free.
    This is clear as any homomorphic image of $K_{r+1}$ is $K_{r+1}$.
\end{proof}
\begin{remark}
    In the previous section, we showed that \cref{thm:entropic-turan} implies the density Tur\'an theorem using a simpler argument.
    This turns out to be the direction we care about in this paper.
    For all the Tur\'an-type results proven later in this paper using entropy and \cref{prop:entropic-lagrangian}, we may also avoid the use of \cref{prop:entropic-lagrangian} by a similar simpler argument.
    However, we think \cref{prop:entropic-lagrangian} is interesting on its own, so we establish the proposition here and will freely use it from now on.
\end{remark}

    Setting $p=2$, we can now prove \cref{thm:spectral-Turan-tree} by applying \cref{thm:entropic-turan} and sampling a random homomorphic image of the tree $T$ in a way similar to \cref{claim:StarSido}.
    \begin{proof}[Proof of \cref{thm:spectral-Turan-tree}]
        From \cref{prop:entropic-lagrangian} and the observation in \cref{ex:entropic-lagrangian}, there exists a random edge with uniform ordering $(X,Y)$ on $G$ such that $\log_2\rho(G)=\HH(Y\mid X)$. By \cref{thm:entropic-turan}, we have
        \[\ell\log_2\rho(G)=\ell\HH(Y\mid X)\leq \HH(X)+(\ell-1)\HH(Y\mid X)+\log_2\left(1-\frac{1}{r}\right).\]
        
        Let $v_1,\ldots, v_{\ell}$ be an ordering of the vertices of $T$ where for every $i\in\{2,\ldots,\ell\}$,  the vertex $v_i$ is adjacent to exactly one $v_j$ with $j<i$.
        Now, we sample random vertices $X_1,\dots,X_\ell$ in $G$ as follows. Let $X_1$ be a random vertex sampled using the law of $X$. Assume we have already sampled $X_1,\dots,X_{i-1}$, and assume $v_j$ is the neighbor of $v_i$ with $j<i$. We sample $X_i$ conditionally independently such that $X_i\mid X_j\sim Y\mid X$. It follows that $X_1,\dots,X_\ell$ is always a homomorphic image of $T$ in $G$. Also, from the way we sample, we know that $\HH(X_1,\dots,X_\ell)=\HH(X)+(\ell-1)\HH(Y\mid X)$. Thus, we have
        \[\HH(X)+(\ell-1)\HH(Y\mid X)=\HH(X_1,\dots,X_\ell)\leq \log_2\#\{\text{homomorphisms from $T$ to $G$}\},\]
        and we are done by combining this with the previous inequality and rearranging.
    \end{proof}

    For general $p$, recall that our definition of $b_p(G)$ matches the definition of $p$-spectral radius given by Keevash, Lenz and Mubayi. Thus, by combining \cref{prop:entropic-lagrangian} with \cref{thm:entropic-turan}, we recover the following theorem for graphs by Kang and Nikiforov \cite{KN14}.
    \begin{theorem}[\cite{KN14}]\label{thm:p-spectral-Turan}
        Let $r\geq 2$ be a positive integer and $p\geq 1$ be a real number. For any $K_{r+1}$-free graph $G$ with $n$ vertices and $m$ edges, we have
        \[b_p(G)\leq \left(1-\frac{1}{r}\right)n^{2-2/p},\]
        and
        \[b_p(G)\leq \left(1-\frac{1}{r}\right)^{1/p}(2m)^{1-1/p}.\]
    \end{theorem}
    \begin{proof}
        From \cref{prop:entropic-lagrangian}, there exists a random edge with uniform ordering $(X,Y)$ on $G$ such that $\log_2b_p(G)=\HH(X,Y)-\frac{2}{p}\HH(X)$. We have
        \[\HH(X,Y)-\frac{2}{p}\HH(X)\leq \left(2-\frac{2}{p}\right)\HH(X)+\log_2\left(1-\frac{1}{r}\right)\leq \left(2-\frac{2}{p}\right)\log_2 n+\log_2\left(1-\frac{1}{r}\right),\]
        and
        \[\HH(X,Y)-\frac{2}{p}\HH(X)\leq \left(1-\frac{1}{p}\right)\HH(X,Y)+\frac{1}{p}\log_2\left(1-\frac{1}{r}\right)\leq \left(1-\frac{1}{p}\right)\log_2 (2m)+\frac{1}{p}\log_2\left(1-\frac{1}{r}\right).\qedhere\]
    \end{proof}
    
    We also remark that, by utilizing \cref{prop:entropic-lagrangian}, we can translate \cref{thm:main-entropy-ver} and also results in \cref{sec:known} into spectral results using arguments in the proofs of \cref{thm:spectral-Turan-tree} and \cref{thm:p-spectral-Turan}.

\section{Partial hypergraphs}\label{sec:partial-hypergraph}
In this section, we introduce some notations and an entropic lemma that will be useful in the later sections.
Those notations are non-standard and are set for our own notational convenience when describing hypergraphs and homomorphisms.

A \emph{partial $k$-graph} $F$ is a simplicial complex whose faces have size at most $k$.
Its set of vertices is denoted by $V(F)$, and its set of \emph{faces}, or \emph{partial edges}, is denoted by $E(F)$.
A \emph{homomorphism} from a partial $k$-graph $F$ to a $k$-graph $G$ is a map $f:V(F)\to V(G)$ such that for any partial edge $e\in E(F)$, $f$ is injective on $e$ and $f(e)$ is contained in some edge in $E(G)$.
Now for any partial $k$-graph $F$, its \emph{extension} $\Tilde{F}$ is the $k$-graph obtained as follows: first let $E'$ be the set of maximal partial edges in $E(F)$, and then extend each partial edge in $E'$ to a $k$-edge by adding in extra vertices, where two different edges do not share any extra vertices.
Notice that if $F$ is a simplicial complex generated by edges of some $k'$-graph $F'$ with $k'<k$, then $\Tilde{F}$ is the extension of $F'$ as defined in the introduction.

\begin{example}[Definition of partial tents]
    In \cref{sec:main-proof}, the partial $k$-graphs and the corresponding extensions of concern would be the following.
For any partition $\lambda$ of $k$ with $\ell\eqdef \ell(\lambda)\geq 2$, the \emph{partial $\lambda$-tent} $\Delta^p_\lambda$ is the partial $k$-graph obtained by taking the simplicial complex generated by $\Delta_\lambda$, and then restricting it to $e\cup \{v\}$ where $e$ is the base and $v$ is the apex.
It is easy to verify that $\Delta_\lambda$ is the extension of the partial $k$-graph $\Delta^p_{\lambda}$.
\end{example}

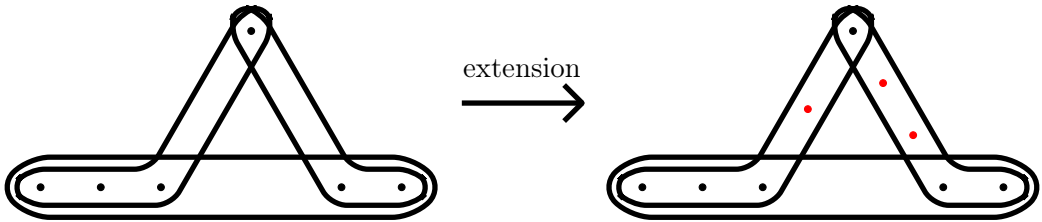
\begin{figure}[h]\centering\label{fig:PartialTent}
\begin{tikzpicture}[scale=0.8]

    \coordinate (A) at (0,0);
    \coordinate (B) at (1,0);
    \coordinate (C) at (2,0);
    \coordinate (D) at (5,0);
    \coordinate (E) at (6,0);
    \coordinate (F) at (4.5,1.732/2);
    \coordinate (G) at (2.75,1.732*3/4);
    \coordinate (H) at (4,1.732);
    \coordinate (I) at (3.5,1.732*3/2);

    \draw [fill] (A) circle (1.6pt);
    \draw [fill] (B) circle (1.6pt);
    \draw [fill] (C) circle (1.6pt);
    \draw [fill] (D) circle (1.6pt);
    \draw [fill] (E) circle (1.6pt);
    \draw [fill] (I) circle (1.6pt);
% Define a path with rounded corners
\draw[rounded corners=8pt,black,line width=2pt] 
    (0-0.2,0.5)--(-0.5-0.2,0)--(0-0.2,-0.5)--(6+0.2,-0.5)--(6.5+0.2,0)--(6+0.2,0.5)--cycle;
    
\draw[rounded corners=6pt,black,line width=2pt] 
(0-0.2,0.3)--(-0.3-0.2,0)--(0-0.2,-0.3)--(2+0.3/1.732,-0.3)--(3.5+0.1+0.15*1.732,1.732*3/2+0.1*1.732-0.15)--(3.5+0.1+0.15,1.732*3/2+0.1*1.732+0.15*1.732)--(3.5+0.1-0.15*1.732,1.732*3/2+0.1*1.732+0.15)--(2-0.3/1.732,0.3)--cycle;

\draw[rounded corners=6pt,black,line width=2pt] 
(6+0.2,0.3)--(6+0.3+0.2,0)--(6+0.2,-0.3)--(5-0.3/1.732,-0.3)--(3.5-0.1-0.15*1.732,1.732*3/2+0.1*1.732-0.15)--(3.5-0.1-0.15,1.732*3/2+0.1*1.732+0.15*1.732)--(3.5-0.1+0.15*1.732,1.732*3/2+0.1*1.732+0.15)--(5+0.3/1.732,0.3)--cycle;

\draw[black,line width=2pt]
(7,1.4)--(9,1.4);

\draw[black,line width=2pt]
(8.7,1.7)--(9,1.4);

\draw[black,line width=2pt]
(8.7,1.1)--(9,1.4);

\node at (8,2) {extension};

    \draw [fill] (9,1.4) circle (0.9pt);

    \draw [fill] (A)+(10,0) circle (1.6pt);
    \draw [fill] (B)+(10,0) circle (1.6pt);
    \draw [fill] (C)+(10,0) circle (1.6pt);
    \draw [fill] (D)+(10,0) circle (1.6pt);
    \draw [fill] (E)+(10,0) circle (1.6pt);
    \draw [fill,red] (F)+(10,0) circle (1.6pt);
    \draw [fill,red] (G)+(10,0) circle (1.6pt);
    \draw [fill,red] (H)+(10,0) circle (1.6pt);
    \draw [fill] (I)+(10,0) circle (1.6pt);
% Define a path with rounded corners
\draw[rounded corners=8pt,black,line width=2pt] 
    (10+0-0.2,0.5)--(10-0.5-0.2,0)--(10+0-0.2,-0.5)--(10+6+0.2,-0.5)--(10+6.5+0.2,0)--(10+6+0.2,0.5)--cycle;
    
\draw[rounded corners=6pt,black,line width=2pt] 
(10+0-0.2,0.3)--(10-0.3-0.2,0)--(10+0-0.2,-0.3)--(10+2+0.3/1.732,-0.3)--(10+3.5+0.1+0.15*1.732,1.732*3/2+0.1*1.732-0.15)--(10+3.5+0.1+0.15,1.732*3/2+0.1*1.732+0.15*1.732)--(10+3.5+0.1-0.15*1.732,1.732*3/2+0.1*1.732+0.15)--(10+2-0.3/1.732,0.3)--cycle;

\draw[rounded corners=6pt,black,line width=2pt] 
(10+6+0.2,0.3)--(10+6+0.3+0.2,0)--(10+6+0.2,-0.3)--(10+5-0.3/1.732,-0.3)--(10+3.5-0.1-0.15*1.732,1.732*3/2+0.1*1.732-0.15)--(10+3.5-0.1-0.15,1.732*3/2+0.1*1.732+0.15*1.732)--(10+3.5-0.1+0.15*1.732,1.732*3/2+0.1*1.732+0.15)--(10+5+0.3/1.732,0.3)--cycle;

%\node at (3.5,0) {Base};
%\node at (4.7,1.732*3/2) {Apex};
\end{tikzpicture}
\caption{Partial $(3,2)$-tent and its extension. Note that for the partial tent, only the maximal edges are shown.}
\end{figure}

Those definitions are useful as for any partial $k$-graph $F$, a homomorphism $F\to G$ is essentially the same as a homomorphism $\Tilde{F}\to G$.
This would be helpful later as instead of considering homomorphisms from $\Delta_{\lambda}$, we can consider homomorphisms from $\Delta^p_{\lambda}$, which are easier to describe.

\begin{proposition}\label{prop:partial-tent-hom}
    Let $F$ be a partial $k$-graph, and let $G$ be a $k$-graph.
    Then there is a homomorphism from $F$ to $G$ if and only if there is a homomorphism from $\Tilde{F}$ to $G$.
\end{proposition}
\begin{proof}
    For any homomorphism $f:V(\Tilde{F})\to V(G)$ from $\Tilde{F}$ to $G$, its restriction $f|_{V(F)}$ is a homomorphism from $F$ to $G$.
    Conversely, suppose that $g:V(F)\to V(G)$ is a homomorphism from $F$ to $G$.
    Note that for every $e\in E(\Tilde{F})$, we have that $g$ is injective on $e\cap V(F)$ and $g(e\cap V(F))$ is contained in some edge in $G$.
    As any vertex in $V(\Tilde{F})\backslash V(F)$ is in exactly one edge in $E(\Tilde{F})$, it is possible to extend $g$ to $\Tilde{g}:V(\Tilde{F})\to V(G)$ so that $g(e)$ is an edge in $G$ for each $e\in E(\Tilde{F})$.
    The extended map $\Tilde{g}$ is indeed a homomorphism from $\Tilde{F}$ to $G$.
\end{proof}

Later on, as in the proof in \cref{sec:entropy}, we will need to show that we can sample random homomorphisms from some tree-like structures with high entropy.
Before we can do so, we need to first describe what the tree-like structures are.

\begin{definition}[Partial forest and forest sequence]
    
For any partial $k$-graph $F$, any linear order $<$ on $V(F)$, and any vertex $v\in V(F)$, let $M_{F,<}(v)$ be the set of partial edges whose maximum vertex is $v$.
A partial $k$-graph $F$ is a \emph{partial forest} with respect to a linear order $<$ on $V(F)$ if for every $v\in V(F)$, there is exactly one maximal partial edge $e_v$ in $M_{F,<}(v)$.
In this case, the \emph{forest sequence} of $(F,<)$ is a sequence $(n_{1},\ldots, n_{k})$ where for each $i\in[k]$, $n_{i}$ is the number of vertices $v\in V(F)$ with $\abs{e_v}=i$. 
\end{definition}

We also define quantities that are analogs of the quantity $2^{\HH(Y\mid X)-\HH(X)}$ we used in \cref{sec:entropy}.

\begin{definition}[Ratio sequence]
    Let $(X_1,\dots,X_k)\in V(G)^k$ be a random edge with uniform ordering on a $k$-graph $G$. We define the \emph{ratio sequence} $0< x_1\leq\dots\leq x_k= 1$ of $(X_1,\dots,X_k)$ by $x_i=2^{\HH(X_i\mid X_{i+1},\dots,X_k)-\HH(X_i)}$ for each $i\in [k]$.
\end{definition}

We are now ready to sample homomorphisms from partial forests with high entropy.

\begin{lemma}\label{lemma:sample-tree}
    Let $(X_1,\ldots, X_k)$ be a random edge with uniform ordering on a $k$-graph $G$ and let $x_1,\ldots, x_k$ be its ratio sequence.
    For any partial forest $F$ with a linear order $<$, if $(n_1,\ldots, n_k)$ is its forest sequence, then one can sample a random homomorphism $(Y_v)_{v\in V(F)}$ from $F$ to $G$ with entropy equal to
    \[v(F)H(X_1)+\log_2\left(\prod_{i=1}^{k}x_i^{n_{k+1-i}}\right).\]
    Moreover, the random homomorphism can be sampled such that for any partial edge $e\in E(F)$, the distribution of $(Y_v)_{v\in e}$ is the same as $(X_i)_{k-\abs{e}+1\leq i\leq k}$.
\end{lemma}
\begin{proof}
    We will induct on $v(F)$.
    The case $v(F)=0$ is vacuously true.
    Now suppose that it holds for partial forest of size $v(F)-1$.
    Let $v_{\max}$ be the maximum vertex in $V(F)$.
    Then $F\backslash \{v_{\max}\}$ is also a partial forest, and so we may sample a random homomorphism $(Y_v)_{v\in V(F)\backslash \{v_{\max}\}}$ with the prescribed properties.
    Let $e$ be the maximal partial edge in $M_{F,<}(v_{\textup{max}})$, and let $j=k+1-\abs{e}$.
    By the inductive hypothesis, $(Y_v)_{v\in e\backslash v_{\max}}$ is identically distributed as $(X_i)_{j+1\leq i\leq k}$.
    Therefore, we may sample $Y_{v_{\max}}$ given $(Y_v)_{v\in e\backslash v_{\max}}$ conditionally independently so that $(Y_v)_{v\in e}$ is identically distributed as $(X_i)_{j\leq i\leq k}$.
    This way,
    \begin{align*}
        \HH\left((Y_v)_{v\in V(F)}\right)=&\HH\left((Y_v)_{v\in V(F)\backslash \{v_{\max}\}}\right)+\HH\left(Y_{v_{\max}}\mid (Y_v)_{v\in e\backslash \{v_{\max}\}}\right)\\
        =&\left(v(F)-1\right)H(X_1)+\log_2\left(x_j^{-1}\prod_{i=1}^{k}x_i^{n_{k+1-i}}\right)+H(X_j\mid X_{j+1},\ldots, X_k)\\
        =&v(F)H(X_1)+\log_2\left(\prod_{i=1}^{k}x_i^{n_{k+1-i}}\right)
    \end{align*}
    where we use that $\HH(X_i)= \HH(X_1)$ for any $i\in[k]$.
    It remains to show that for any partial edge $e'$ containing $v_{\max}$, the distribution of $(Y_v)_{v\in e'}$ is the same as $(X_i)_{k-\abs{e'}+1\leq i\leq k}$.
    This is true as $e'\subseteq e$ by the definition of $e$ and $v_{\max}$, and the distribution $(X_1,\ldots, X_k)$ is symmetric.
\end{proof}

\section{Proof of \cref{thm:main-tent,thm:one-tent}}\label{sec:main-proof}
In this section, we will first give two proofs of \cref{thm:main-tent}.
We will then show how \cref{thm:main-tent} implies \cref{thm:one-tent}.
Finally, we will conclude this section with a proof of \cref{thm:not-Turan}.

Throughout this section, we will fix a $k$-graph $G$ and a random edge with uniform ordering $(X_1,\ldots, X_k)$ on $G$.
We will also set $0< x_1\leq \cdots\leq  x_k=1$ to be its ratio sequence.
We make an observation that to upper bound $b(G)=b_{\textup{entropy}}(G)$, it suffices to upper bound\linebreak $2^{\HH(X_1,\ldots, X_k)-k\HH(X_1)}=x_1\cdots x_{k-1}$ by the chain rule. 
Therefore, the upper bound of \cref{thm:main-tent} follows from the following statement.

\begin{theorem}\label{thm:main-entropy-ver}
    If $G$ is $\lambda$-tent-hom-free for every $\abs{\lambda} = k$ and $\ell(\lambda)=2$, then we have
    \[\HH(X_1,\dots,X_k)- k\HH(X_1)=\log_2(x_1\cdots x_k)\leq \log_2\frac{k!}{k^k}.\]
\end{theorem}

We first show that \cref{thm:main-tent} indeed follows from \cref{thm:main-entropy-ver}.
\begin{proof}[Proof of \cref{thm:main-tent} using \cref{thm:main-entropy-ver}]
    First, it is clear that $\pi(\cF_k)\geq k!/k^k$ as a single edge does not contain any homomorphic image of any tents, and it has blowup density $k!/k^k$.
    To show the reverse inequality, if $G$ is $\cF_k$-hom-free, then by \cref{thm:main-entropy-ver}, we have $b(G)=b_{\textup{entropy}}(G)\leq k!/k^k$.
    Combining with supersaturation (\cref{thm:supersaturation}), we get $\pi(\cF_k)\leq k!/k^k$.
\end{proof}

\subsection{First proof of \cref{thm:main-entropy-ver}}

To prove \cref{thm:main-entropy-ver}, we will apply \cref{lemma:sample-tree} and \cref{lemma:mix} to obtain several inequalities involving $x_1,\ldots,x_{k}$.
Then we will solve for the maximum of $x_1\cdots x_{k-1}$ subject to the inequalities.

\begin{lemma}\label{lemma:alt-tent-ineq}
If $G$ is $\lambda$-tent-hom-free for every $\abs{\lambda} = k$ and $\ell(\lambda)=2$, then for any $i,j\in[k]$ with $i+j\leq k$, we have $x_i+x_j\leq x_{i+j}$.
    
\end{lemma}
\begin{proof}
    We will consider two partial forests $F^{(1)}$ and $F^{(2)}$ both on $V=\{v_1,\ldots, v_k,w\}$.
    Let $F^{(1)}$ be spanned by the two partial edges $\{v_1,\ldots, v_k\}$ and $\{v_{i+1},\ldots,v_k, w\}$.
    Let $F^{(2)}$ be spanned by the two partial edges $\{v_1,\ldots, v_k\}$ and $\{v_1,\ldots, v_{k-j}, w\}$.
    Then both partial $k$-graphs are indeed partial forests with respect to the linear order $v_1<\cdots <v_k<w$.
    It is clear that in $F^{(1)}$ with the forest sequence $(n_1,\ldots, n_k)$, the vertices $v_1,\ldots, v_k$ contribute one to $n_1,\ldots, n_k$, respectively, and $w$ contributes to $n_{k-i+1}$.
    Similarly, the forest sequence of $F^{(2)}$ is all-one except for $n_{k-j+1}=2$.

    Let $(Y^{(1)}_v)_{v\in V},(Y^{(2)}_v)_{v\in V}$ be the random homomorphism from $F^{(1)},F^{(2)}$ given by \cref{lemma:sample-tree}, respectively.
    Note that if some tuple of vertices is in the supports of both $(Y^{(1)}_v)_{v\in V}$ and $(Y^{(2)}_v)_{v\in V}$, then this tuple corresponds to a homomorphism from $F^{(1)}\cup F^{(2)}$ to $G$.
    As $F^{(1)}\cup F^{(2)}$ clearly contains a partial $(i,k-i)$-tent with base $\{v_1,\ldots, v_k\}$ and apex $w$, we know that the two random homomorphisms have disjoint support.
    Suppose that $(Z_v)_{v\in V}$ is the mixture given by \cref{lemma:mix}, then by \cref{lemma:mix,lemma:sample-tree} we know
    \[\left(x_1\cdots x_{k-1}\cdot x_i+x_1\cdots x_{k-1}\cdot x_j\right)2^{(k+1)\HH(X_1)}\leq 2^{\HH((Z_v)_{v\in V})}.\]
    Observe that both $F^{(1)}$ and $F^{(2)}$ contains the partial edges $\{v_1,\ldots, v_k\}$ and $\{v_{i+1},\ldots, v_{k-j},w\}$.
    Therefore $(Y^{(1)}_{v_1},\ldots, Y^{(1)}_{v_k})$ and $(Y^{(2)}_{v_1},\ldots, Y^{(2)}_{v_k})$ both have the same distributions as $(X_1,\ldots, X_k)$ by \cref{lemma:sample-tree}, which shows that $(Z_{v_1},\ldots, Z_{v_k})$ has the same distribution as $(X_1,\ldots, X_k)$ as well.
    Using a similar argument, we can show that $(Z_w, Z_{v_{i+1}},\ldots, Z_{v_{k-j}})$ has the same distribution as $(X_{i+j}, \ldots, X_k)$.
    As a consequence,
    \begin{align*}
        \HH\left((Z_v)_{v\in V}\right)\leq& \HH(Z_{v_1},\ldots, Z_{v_k})+\HH(Z_w\mid Z_{v_{i+1}},\ldots, Z_{v_{k-j}})\\
        =&\HH(X_1,\ldots, X_{k})+\HH(X_{i+j}\mid X_{i+j+1},\ldots, X_{k})\\
        =& (k+1)\HH(X_1)+\log_2(x_1\cdots x_{k-1}\cdot x_{i+j}).
    \end{align*}
    This shows that 
    \[x_1\cdots x_{k-1}2^{(k+1)\HH(X_1)}(x_i+x_j)\leq x_1\cdots x_{k-1}2^{(k+1)\HH(X_1)}\cdot x_{i+j}\]
    and so the desired statement follows.
\end{proof}

Our next goal is to upper bound $ x_1\cdots x_{k-1}$ and show that the maximum is obtained when $x_i = i/k$ for every $i=1,\ldots, k-1$.
To upper bound the product, we prove the following auxiliary inequality.

\begin{lemma}\label{lemma:aux-ineq}
    Suppose that $y_1,\ldots, y_k$ are some non-negative real numbers with $y_i+y_j\leq y_{i+j}$ for any $i,j\in[k]$ with $i+j\leq k$.
    Then
    \[y_1\cdots y_k\leq k!\left(\frac{y_1+\cdots +y_k}{\binom{k+1}{2}}\right)^k.\]
\end{lemma}
\begin{proof}
    We will prove this by induction.
    It clearly holds when $k=1$.
    Now suppose that $k\geq 2$ and the statement holds for $k-1$.
    Then by the inductive hypothesis,
    \[y_1\cdots y_k\leq (k-1)!\left(\frac{y_1+\cdots +y_{k-1}}{\binom{k}{2}}\right)^{k-1}y_k\leq k!\left(\frac{(k-1)\cdot\frac{y_1+\cdots+y_{k-1}}{\binom{k}{2}}+\frac{y_k}{k}}{k}\right)^k\]
    by AM-GM.
    Since
    \[y_1+\dots+y_{k-1} = \frac{1}{2}\sum_{i=1}^{k-1}(y_i+y_{k-i})\leq \frac{k-1}{2}y_k,\]
    we know 
    \[(k-1)\cdot\frac{y_1+\cdots+y_{k-1}}{\binom{k}{2}}+\frac{y_k}{k} = \frac{2}{k}\left(y_1+\cdots +y_{k-1}+\frac{y_k}{2}\right)\leq \frac{2}{k}\cdot \frac{k}{k+1}\left(y_1+\cdots +y_k\right)\]
    and so
    \[y_1\cdots y_k\leq k!\left(\frac{\frac{2}{k+1}(y_1+\cdots+y_k)}{k}\right)^k =k!\left(\frac{y_1+\cdots +y_k}{\binom{k+1}{2}}\right)^k, \]
    as desired.
\end{proof}

Combining \cref{lemma:alt-tent-ineq} and \cref{lemma:aux-ineq}, we are now ready to prove \cref{thm:main-entropy-ver}.

\begin{proof}[Proof of \cref{thm:main-entropy-ver}]
    By \cref{lemma:alt-tent-ineq}, $x_1,\ldots, x_k$ are non-negative reals satisfying the condition of \cref{lemma:aux-ineq}.
    We also know that $x_k=1$, so $x_1+\cdots +x_k\leq \frac{k-1}{2}+1 = \frac{k+1}{2}$.
    Thus by \cref{lemma:aux-ineq},
    \[x_1\cdots x_{k-1} = x_1\cdots x_k\leq k!\left(\frac{\frac{k+1}{2}}{\binom{k+1}{2}}\right)^k = \frac{k!}{k^k},\]
    which is the desired statement
\end{proof}
\subsection{Second proof of \cref{thm:main-entropy-ver}}
Here, we give an alternative proof using much more complicated partial forests.
Although the proof is more involved, this proof would be the one we generalize later in \cref{sec:known}.

\begin{lemma}\label{lemma:tent-ineq}
    If $G$ is $\lambda$-tent-hom-free for every $\abs{\lambda}=k$ and $\ell(\lambda)=2$, then for every $i\in[k-1]$, we have $x_j<x_{i+1}$ for each $j\leq i$ and
    \[\prod_{j=1}^{i}\frac{x_j}{x_{i+1}-x_j}\leq 1.\]
\end{lemma}
\begin{proof}

    We will fix $i$ throughout this proof.
    As in what we did in \cref{sec:entropy}, we will temporarily fix an integer $N\in\NN$ that will later be taken to infinity.
    For any $1=t_0<t_1<t_2<\cdots <t_{i+1}= N+1$, we will define a partial forest $F^{(\vec{t})}$ on $V=\{v_1,\ldots, v_{k-i-1}, w_1,\ldots, w_N\}.$
    The partial forest $F^{(\vec{t})}$ is spanned by the partial edges $\{v_1,\ldots, v_{k-i-1},w_m, w_{t_{j+1}},\ldots, w_{t_i}\}$ for every $t_j\leq m<t_{j+1}$.
    This is indeed a partial forest with respect to the linear order $<$ with $v_1<\cdots<v_{k-i-1}<w_N<\cdots <w_1$.
    We can compute the forest sequence with respect to the linear order as follows: each $v_j$ contributes one to $n_j$ for each $j\leq k-i-1$, and each $w_m$ with $t_j\leq m<t_{j+1}$ contributes $1$ to $n_{k-j}$.
    Therefore the forest sequence $(n_1,\ldots, n_k)$ is $(1,\dots,1,t_{i+1}-t_i,\dots,t_1-t_0)$.
    Now let $(Y^{(\vec{t})}_v)_{v\in V}$ be the random homomorphism produced by \cref{lemma:sample-tree}.
    This gives
    \begin{align}\label{eq:tent-2ndproof-eq1}
        \HH\left((Y^{(\vec{t})}_v)_{v\in V}\right) = (N+k-i-1)\HH(X_1)+\log_2\left(x_{i+2}\cdots x_{k}\cdot \prod_{j\leq i+1}x_j^{t_{j}-t_{j-1}}\right).
    \end{align}

    We will now show that the supports of $(Y^{(\vec{t})}_v)_{v\in V}$ are disjoint for different choices of $\vec{t}$.
    Suppose for the sake of contradiction that for some $\vec{t}\neq \vec{t}'$ there is a tuple of vertices from $V(G)$ lying in the supports of $(Y^{(\vec{t})}_v)_{v\in V}$ and $(Y^{(\vec{t}')}_v)_{v\in V}.$
    Then this tuple witnesses a homomorphism sending $F^{(\vec{t})}\cup F^{(\vec{t}')}$ to $G$.
    We will show a contradiction by demonstrating that $F^{(\vec{t})}\cup F^{(\vec{t}')}$ contains a homomorphic image of some partial $\lambda$-tent with $\ell(\lambda)=2$.

    Let $j\geq 1$ be the minimum index in which $\vec{t}$ and $\vec{t}'$ differ, and without loss of generality, suppose that $t_j'<t_j$.
    Then we can find partial edges $e=\{v_1,\ldots,v_{k-i-1},w_{t_0},w_{t_1},\ldots, w_{t_i}\}$,\linebreak $e_1=\{v_1,\ldots, v_{k-i-1},w_{t_j'}, w_{t_j},\ldots, w_{t_i}\}$ in $F^{(\vec{t})}$ and $e_2 = \{w_{t_0'},\ldots, w_{t_j'}\}$ in $F^{( \vec{t}')}$.
    By the minimality of $j$, we know $e_2=\{w_{t_0},\ldots, w_{t_{j-1}}, w_{t_j'}\}.$
    Note that $e,e_1,e_2$ form a partial $(k-j,j)$-tent with base $e$ and apex $w_{t_j'}$, showing that $F^{(\vec{t})}\cup F^{(\vec{t}')}$ contains a partial $(k-j,j)$-tent, which is a contradiction.
    
    Therefore we may now apply \cref{lemma:mix} with $a=1$.
    Suppose that $(Z_v)_{v\in V}$ is the resulting mixture of $(Y_v^{(\vec{t})})_{v\in V}$ for all possible $\vec{t}$.
    Then
    \begin{align}\label{eq:tent-2ndproof-eq3}\frac{\sum_{1=t_0<t_1<\cdots<t_{i+1}=N+1}2^{\HH\left((Y_v^{(\vec{t})})_{v\in V}\right)}}{2^{\HH\left((Z_v)_{v\in V}\right)}}\leq 1.\end{align}
    On the other hand, by \cref{lemma:sample-tree} and the fact that $\{v_1,\ldots, v_{k-i-1},w_m\}$ is present in all partial forests we take for any $m\in[N]$, we know that $(Z_{v_1},\ldots, Z_{v_{k-i-1}},Z_{w_m})$ has the same distribution as $(X_{i+1},\ldots, X_k)$ for each $m\in[N]$.
    Hence
    \begin{align}
        \HH\left((Z_v)_{v\in V}\right)\leq& \HH(Z_{v_1},\ldots, Z_{v_{k-i-1}})+\sum_{m=1}^{N}\HH(Z_{w_m}\mid Z_{v_1},\ldots, Z_{v_{k-i-1}}) \nonumber\\
        =& \HH(X_{i+2},\ldots,X_k)+N\HH(X_{i+1}\mid X_{i+2},\ldots,X_k) \nonumber\\
        =& (N+k-i-1)\HH(X_1)+\log_2(x_{i+2}\cdots x_k\cdot x_{i+1}^N).\label{eq:tent-2ndproof-eq2}
    \end{align}
    Thus \cref{eq:tent-2ndproof-eq1,eq:tent-2ndproof-eq2} now gives
    \[\frac{\sum_{\vec{t}}2^{\HH\left((Y_v^{(\vec{t})})_{v\in V}\right)}}{2^{\HH\left((Z_v)_{v\in V}\right)}}\geq \frac{2^{(N+k-i-1)\HH(X_1)}\cdot x_{i+2}\cdots x_k\cdot\sum_{\vec{t}}  \prod_{j\leq i+1}x_j^{t_{j}-t_{j-1}}}{2^{(N+k-i-1)\HH(X_1)}\cdot x_{i+2}\cdots  x_k\cdot x_{i+1}^{N}}=\sum_{\vec{t}}\prod_{j\leq i+1}\left(\frac{x_j}{x_{i+1}}\right)^{t_j-t_{j-1}},\]
    where $\sum_{\vec{t}}$ sums over $\vec{t}$ with $1=t_0<t_1<\cdots<t_{i+1}=N+1$.
    Note that we may replace $j\leq i+1$ by $j\leq i$ in the product.
    If we take $\delta_j = t_j-t_{j-1}\in\NN$, then
    \[\sum_{\vec{t}}\prod_{j\leq i}\left(\frac{x_j}{x_{i+1}}\right)^{t_j-t_{j-1}}=\sum_{\delta_1+\dots+\delta_i\leq N}\prod_{j\leq i}\left(\frac{x_j}{x_{i+1}}\right)^{\delta_j}.\]
    Since the sum is always upper bounded by $1$ when we take $N$ goes to infinity, we must have $x_j<x_{i+1}$ for each $j\in[i]$. 
    Thus, we have
    \[\lim_{N\to\infty}\sum_{\vec{t}}\prod_{j\leq i}\left(\frac{x_j}{x_{i+1}}\right)^{t_j-t_{j-1}}=\sum_{\delta_1,\ldots, \delta_{i}\in \NN}\prod_{j\leq i}\left(\frac{x_j}{x_{i+1}}\right)^{\delta_j} =\prod_{j\leq i}\frac{x_j}{x_{i+1}-x_j}.\]
    Therefore
    \[\liminf_{N\to\infty}\frac{\sum_{1=t_0<t_1<\cdots<t_{i+1}=N+1}2^{\HH\left((Y_v^{(\vec{t})})_{v\in V}\right)}}{2^{\HH\left((Z_v)_{v\in V}\right)}}\geq \prod_{j\leq i}\frac{x_j}{x_{i+1}-x_j}.\]
    Combining this with \cref{eq:tent-2ndproof-eq3} gives the desired inequality.
\end{proof}

Once again, to prove \cref{thm:main-entropy-ver}, we need to upper bound $x_1\cdots x_{k-1}$ given the inequalities in \cref{lemma:tent-ineq} and show that it is maximized when $x_i = i/k$ for every $i=1,\ldots, k-1$.
We will prove a slightly stronger statement, which will also be useful in the next section.

\begin{lemma}\label{lemma:another-aux-ineq}
    Let $k$ be a positive integer.
    Fix real numbers $0<z_1<\cdots <z_k$.
    Let $0< y_1< \ldots<  y_k$ be real numbers with
    \[\prod_{j\leq i}\frac{y_j}{y_{i+1}-y_j}\leq \prod_{j\leq i}\frac{z_j}{z_{i+1}-z_j}\]
    for any $i=1,\ldots, k-1$.
    Then 
    \[y_1\cdots y_{k-1}\leq \frac{z_1\cdots z_{k-1}}{z_k^{k-1}}y_k^{k-1}.\]
\end{lemma}
\begin{proof}
    We will prove by induction on $k$.
    When $k=1$ this is clearly true.
    Now suppose that $k\geq 2$ and the statement is true for all smaller $k$.
    Then we have
    \[\frac{y_1\cdots y_i}{z_1\cdots z_i}\leq \frac{y_{i+1}^i}{z_{i+1}^{i}}\]
    for all $i<k-1$ by the inductive hypothesis.
    Now let 
    \[\alpha_i =\frac{1}{i}\sum_{j\leq i} \frac{z_j}{z_k-z_j}\] for any $i\leq k-1$.
    Note that for any $i<k-1$, we have
    \[\left(\frac{y_1\cdots y_{i+1}}{z_1\cdots z_{i+1}}\right)^{\alpha_{i+1}}\leq \left(\frac{y_1\cdots y_i}{z_1\cdots z_i}\right)^{\alpha_{i}}\left(\frac{y_{i+1}^i}{z_{i+1}^i}\right)^{(\alpha_{i+1}-\alpha_i)}\left(\frac{y_{i+1}}{z_{i+1}}\right)^{\alpha_{i+1}}=\left(\frac{y_1\cdots y_i}{z_1\cdots z_i}\right)^{\alpha_{i}}\left(\frac{y_{i+1}}{z_{i+1}}\right)^{\frac{z_{i+1}}{z_k-z_{i+1}}}.\]
    Here, we are using that $\alpha_{i+1}-\alpha_i\geq 0$ as $\frac{z_1}{z_k-z_1}<\cdots <\frac{z_{i+1}}{z_k-z_{i+1}}$.
    Multiplying these up for $i=1,\ldots, k-2$, and we get
    \[\left(\frac{y_1\cdots y_{k-1}}{z_1\cdots z_{k-1}}\right)^{\alpha_{k-1}}\leq \left(\frac{y_1}{z_1}\right)^{\frac{z_1}{z_k-z_1}}\cdots\left(\frac{y_{k-1}}{z_{k-1}}\right)^{\frac{z_{k-1}}{z_k-z_{k-1}}}.\]
    Thus
    \begin{align*}
        \left(\frac{y_1\cdots y_{k-1}}{z_1\cdots z_{k-1}}\right)^{\alpha_{k-1}+1}\leq& \prod_{i=1}^{k-1}\left(\frac{y_i}{z_i}\right)^{\frac{z_i}{z_k-z_i}}\left(\frac{y_k-y_i}{z_k-z_i}\right)\\
        =&\prod_{i=1}^{k-1}\left[\left(\frac{y_i}{z_i}\right)^{\frac{z_i}{z_k}}\left(\frac{y_k-y_i}{z_k-z_i}\right)^{\frac{z_k-z_i}{z_k}}\right]^{\frac{z_k}{z_k-z_i}}\\
        \leq &\prod_{i=1}^{k-1}\left(\frac{z_i}{z_k}\cdot \frac{y_i}{z_i}+\frac{z_k-z_i}{z_k}\cdot \frac{y_k-y_i}{z_k-z_i}\right)^{\frac{z_k}{z_k-z_i}}&\textup{(weighted AM-GM)}\\
        = &\prod_{i=1}^{k-1}\left(\frac{y_k}{z_k}\right)^{\frac{z_k}{z_{k}-z_i}}
        =\left(\frac{y_k}{z_k}\right)^{(k-1)(\alpha_{k-1}+1)},
    \end{align*}
    completing the inductive step.
    
\end{proof}

\begin{proof}[Alternative Proof of \cref{thm:main-entropy-ver}]
    Suppose $G$ is $\cF_k$-hom-free.
    Set $z_i=i$ for each $i\in[k]$.
    By \cref{lemma:tent-ineq}, we know that
    \[\prod_{j\leq i}\frac{x_j}{x_{i+1}-x_j}\leq 1 = \prod_{j\leq i}\frac{j}{(i+1)-j}=\prod_{j\leq i}\frac{z_j}{z_{i+1}-z_j}.\]
    Therefore by \cref{lemma:another-aux-ineq} and the fact that $x_k=1$,
    \[x_1\cdots x_{k-1}\leq \frac{(k-1)!}{k^{k-1}}=\frac{k!}{k^k},\]
    as desired.
\end{proof}

\subsection{Proof of \cref{thm:one-tent,thm:not-Turan}}
As mentioned in the introduction, \cref{thm:one-tent} is an immediate corollary of \cref{thm:main-tent}.
We give a detailed argument of how \cref{thm:one-tent} follows from \cref{thm:main-tent} below.

\begin{proof}[Proof of \cref{thm:one-tent}]
    Let $\lambda$ be a partition of $k$ with $\lambda_1\leq \lceil k/2\rceil$ and $\lambda_i = 1$ for all $1<i\leq \ell(\lambda)$.
    Again, it is clear that $\pi(\Delta_\lambda)\geq k!/k^k$, so it suffices to show that $\pi(\Delta_\lambda)\leq k!/k^k$.
    By \cref{thm:supersaturation}, it suffices to show that any $\Delta_\lambda$-hom-free $k$-graph $G$ is also $\Delta_{\lambda'}$-hom-free for any $\lambda'$ with $\abs{\lambda'}=k$ and $\ell(\lambda')=2$.
    This will follow immediately if we show that $\Delta_{\lambda}$ admits a homomorphism to $\Delta_{\lambda'}$ for any such $\lambda'$.
    By \cref{prop:partial-tent-hom}, it is sufficient to show that $\Delta^p_{\lambda}$  admits a homomorphism to $\Delta_{\lambda'}$ for any $\lambda'$ with $\abs{\lambda'}=k$ and $\ell(\lambda')=2$.
    This is now simple: suppose that $\Delta_{\lambda'}$ has base $e'$ and apex $v'$, and $e_1', e_2'$ are two edges such that $\abs{e_i'\cap e'}=\lambda_i'$ for $i\in[2]$.
    We also suppose that $\Delta^p_{\lambda}$ has base $e$ and apex $v$, and $e_1,\ldots, e_{\ell}$ are partial edges such that $\abs{e_i\cap e}=\lambda_i$ for $i\in[\ell]$.
    As $\lambda_1'\geq \lceil k/2\rceil \geq \lambda_1$, we can take $f:e\cup \{v\}\to V(\Delta_{\lambda'})$ so that $f(v)=v'$, $f(e)=e'$ and $f(e\cap e_1)\subseteq e'\cap e_1'$.
    This is a homomorphism from $\Delta^p_{\lambda}$ to $\Delta_{\lambda'}$ as any vertex in $e'$ shares an edge with $v'$ in $\Delta_{\lambda'}$.
\end{proof}

Finally, we give a proof of \cref{thm:not-Turan} by demonstrating a $k$-graph $G$ that has $b(G)>k!/k^k$ and is $\Delta_\lambda$-free for large $\lambda_1$.
Similar to an earlier lower-bound construction by Frankl and F\"uredi \cite{FF89} for $\Delta_{(k-1,1)}$, we will do so by constructing a $k$-graph $G$ so that the intersection of any two edges is small.

\begin{proof}
    Let $\alpha<1$ be some constant that is close to $1$.
    In particular, assume that $\alpha > 1/2$.
    Let $G_{\textup{aux}}$ be an auxiliary graph with vertices $\binom{[2k]}{k}$, and two vertices are connected if the corresponding subsets have intersection at least $\alpha k$.
    Then $G_{\textup{aux}}$ is a regular graph with degree
    \[\sum_{i\leq (1-\alpha)k}\binom{k}{i}^2< k\binom{k}{\lfloor (1-\alpha)k\rfloor}^2 = 2^{(2h(\alpha)+o(1))k},\]
    where $h(\alpha) = -\alpha\log_2\alpha-(1-\alpha)\log_2\alpha$ and we use that 
    \[\binom{k}{(1-\alpha+o(1))k}=2^{(h(\alpha)+o(1))k}\]
    when $\alpha>1/2$.
    
    By the Caro--Wei theorem, there exists an independent set of size
    \[\frac{\binom{2k}{k}}{2^{(2h(\alpha)+o(1))k}} = 2^{(2-2h(\alpha)+o(1))k}.\]
    This corresponds to a $k$-graph $G$ on $[2k]$ with $2^{(2-2h(\alpha)+o(1))k}$ edges so that any two edges have intersection less than $\alpha k$.

    Now if $G$ contains a homomorphic image of $\Delta_\lambda$ where $\lambda_1>\alpha k$, let $e$ be its base and let $e_1$ be the edge with $\abs{e\cap e_1}=\lambda_1$.
    Also let $f$ be a homomorphism from $\Delta_\lambda$ to $G$.
    Then $\abs{f(e)\cap f(e_1)}>\alpha k$, and so $f(e)=f(e_1)$.
    This shows if $v$ is the apex of $\Delta_\lambda$, then $f(v)=f(u)$ for some $u\in e$.
    However, $\{uv\}$ is contained in some edge in $\Delta_\lambda$, which is a contradiction.
    Thus $\pi(\Delta_\lambda)$ is at least $b(G)$, which is at least the density of $G$.
    The density of $G$ is
    \[\frac{k! \cdot 2^{(2-2h(\alpha)+o(1))k} }{(2k)^k} = 2^{(1-2h(\alpha)+o(1))k}\cdot \frac{k!}{k^k},\]
    which is strictly greater than $k!/k^k$ for sufficiently large $k$ as long as $h(\alpha)<1/2$.
    As $h$ is continuous on $[1/2, 1]$ and $h(1)=0$, this is true for $\alpha$ sufficiently close to $1$.
\end{proof}
The proof roughly gives $\alpha \approx 0.89$.
Although our proof is not fully optimized, we believe that it would not give the correct upper bound for $\alpha$ even after being fully optimized.
Therefore we do not pursue this direction.

\section{Other applications of our method}\label{sec:known}

Recall from the introduction that Mubayi \cite{Mub06} showed $\pi(E^{(k)}_{k+1})=k!/k^k$ where $E^{(k)}_{k+1}$ is the extended clique of size $k+1$, and Mubayi and Pikhurko \cite{MP07} strengthened it to $\pi(\Delta_{(1,1,\ldots,1)}) = k!/k^k$.
In fact they both proved more general results than this: Mubayi showed that for each $r\geq k$,
\[\pi(E^{(k)}_{r+1})=b(K^{(k)}_r) = \prod_{i=1}^{k-1}\left(1-\frac{i}{r}\right)\]
and Mubayi and Pikhurko strengthened it as follows: consider the partial $k$-graph $F$ on $r+1$ vertices generated by $[k]$ and all the $2$-subsets of $[r+1]$, and then take its extension $\Tilde{F}$.
Then $\pi(\Tilde{F})=b(K^{(k)}_r)$ as well.
Note that $E^{(k)}_{r+1}$ is the extension of $K_{r+1}$ as a partial $k$-graph, and there is a homomorphism from $K_{r+1}$ to $\Tilde{F}$.
Therefore $\pi(E^{(k)}_{r+1})\leq \pi(\Tilde{F})$ (where we use \cref{thm:supersaturation}), and so $\pi(\Tilde{F}) = b(K^{(k)}_r)$ is indeed a stronger statement.
We remark that Keevash's adaptation \cite[Theorem 3.1]{Kee11} of Sidorenko's argument \cite{Sid89} gives a much more general result than Mubayi and Pikhurko's result in this case, and we refer the readers to Keevash's survey for the statement.

We are able to prove $\pi(\Tilde{F}) = b(K^{(k)}_r)$ as well, though our proof is considerably more complicated, and it seems hard to produce a clean stronger statement.
We nonetheless outline the argument here for readers interested in improving our argument.

\begin{theorem}\label{thm:fix-k-big-r-full-edge}
    Let $k,r$ be positive integers with $r\geq k$.
    For any positive integer $N$, any $i=1,\ldots,k-1$ and any sequence $\vec{t}=(t_0,\ldots, t_{i+1})$ with $1=t_0<t_1<\cdots <t_{i+1}=N+1$, set $F^{(\vec{t})}$ to be the partial forest on $\{v_1,\ldots, v_{k-i-1},w_1,\ldots, w_N\}$ spanned by partial edges $\{v_1,\ldots, v_{k-i-1},w_m, w_{t_{j+1}},\ldots, w_{t_i}\}$ for every $t_j\leq m<t_{j+1}$.
    
    Let $\cF$ be a family of $k$-graphs such that for all positive integers $N$ and $i=1,\ldots, k-1$, if we take the union of any $\binom{r-k+i}{i}+1$ different partial forests $F^{(\vec{t})}$ with $1=t_0<t_1<\cdots <t_{i+1}= N+1$, then its extension is not $\cF$-hom-free.
    
    Then $\pi(\cF) \leq b(K^{(k)}_r)$.
\end{theorem}
\begin{proof}
    Suppose that $G$ is $\cF$-hom-free.
    Let $(X_1,\ldots, X_k)$ be any random edge with uniform ordering on $G$ and let $x_1,\dots,x_k$ be its ratio sequence.
    We first fix some $i\in[k-1]$ and some large positive integer $N$.
    For any $1=t_0<t_1<\cdots <t_{i+1}=N+1$, let $(Y_v^{(\vec{t})})_{v\in V}$ be the random homomorphism from $F^{(\vec{t})}$ to $G$ sampled via \cref{lemma:sample-tree} using the linear order $v_1<\cdots <v_{k-i-1}<w_N<\cdots <w_1$.
    Then by the assumption on $\cF$ and that $G$ is $\cF$-hom-free, we know that the supports of the random homomorphisms $(Y_v^{(\vec{t})})_{v\in V}$ are $\left(\binom{r-k+i}{i}+1\right)$-wise disjoint.
    Therefore, if $(Z_v)_{v\in V}$ is the mixture of the $(Y_v^{(\vec{t})})_{v\in V}$'s provided by \cref{lemma:mix}, we have
    \[ \sum_{1=t_0<t_1<\cdots <t_{i+1}=N+1}2^{\HH\left((Y_v^{(\vec{t})})_{v\in V}\right)}\leq\binom{r-k+i}{i}2^{\HH\left((Z_v)_{v\in V}\right)}.\]

    Recall that from the calculation in the proof of \cref{lemma:tent-ineq}, we have
    \[\liminf_{N\to\infty}\frac{\sum_{1=t_0<t_1<\cdots<t_{i+1}=N+1}2^{\HH\left((Y_v^{(\vec{t})})_{v\in V}\right)}}{2^{\HH\left((Z_v)_{v\in V}\right)}}\geq \prod_{j\leq i}\frac{x_j}{x_{i+1}-x_j}.\]
    Thus, we get
    \[\prod_{j\leq i}\frac{x_j}{x_{i+1}-x_j}\leq \binom{r-k+i}{i}.\]

    Now let $z_i = r-k+i$ for each $i=1,\ldots, k$.
    Then it is easy to verify that
    \[\binom{r-k+i}{i}= \prod_{j\leq i}\frac{z_j}{z_{i+1}-z_j}\]
    for each $i\in[k-1]$.
    Therefore, by \cref{lemma:another-aux-ineq}, we get that
    \[x_1\cdots x_{k-1}\leq \frac{z_1\cdots z_{k-1}}{z_k^{k-1}} = \frac{(r-k+1)\cdots (r-1)}{r^{k-1}} = \prod_{i=1}^{k-1}\left(1-\frac{i}{r}\right)=b(K_r^{(k)}).\]
    This shows that $b(G) = b_{\textup{entropy}}(G)\leq b(K_r^{(k)})$ for any $\cF$-hom-free $k$-graph $G$, and so we have $\pi(\cF)\leq b(K_r^{(k)})$ by \cref{thm:supersaturation}.
\end{proof}

\begin{corollary}\label{cor:fix-k-big-r-full-edge}
    Let $F$ be the partial $k$-graph on $r+1$ vertices generated by $[k]$ and all the $2$-subsets of $[r+1]$.
    Let $\Tilde{F}$ be its extension.
    Then $\pi(\Tilde{F}) = b(K_r^{(k)})$.
\end{corollary}
\begin{proof}
    First of all, it is clear that $K_r^{(k)}$ is $F$-hom-free.
    Therefore, by \cref{prop:partial-tent-hom}, $K_r^{(k)}$ is also $\Tilde{F}$-hom-free, and so $\pi(\Tilde{F})\geq b(K_r^{(k)})$.

    To show that $\pi(\Tilde{F})\leq b(K_r^{(k)})$, it now suffices to show that the assumption of \cref{thm:fix-k-big-r-full-edge} holds for any $i\in[k-1]$.
    Indeed, for any collection $T$ of $\binom{r-k+i}{i}+1$ different possible $\vec{t}$'s, we may construct $S\subseteq\NN$ with size $r-k+i+1$ that satisfies the following: for each $s\in S$ there exists $\vec{t}\in T$ such that $s\in\{t_1,\ldots, t_i\}$, and there exists a $\vec{t}\in T$ with $\{t_1,\ldots, t_i\}\subseteq S$.
    Indeed, set $S' =\bigcup_{\vec{t}\in T}\{t_1,\ldots, t_i\}$.
    Then $\abs{T}\leq \binom{\abs{S'}}{i}$, which shows that $\abs{S'}\geq r-k+i+1$.
    Now simply take $S\subseteq S'$ of size $r-k+i+1$ while containing some $\{t_1,\ldots, t_i\}$ for some $\vec{t}\in T$.
    Label this $\vec{t}$ as $\vec{t^*}$.

    Now we need to show that there is a homomorphic image of $\Tilde{F}$ in the extension of $\bigcup_{\vec{t}\in T}F^{(\vec{t})}$.
    By \cref{prop:partial-tent-hom}, it suffices to construct a homomorphism from $F$ to $\bigcup_{\vec{t}\in T}F^{(\vec{t})}$.
    To do so, we will simply map $1,\ldots,k-i-1$ to $v_1,\ldots, v_{k-i-1}$, map $k-i,\ldots, k$ to $w_{t^*_0},\ldots, w_{t^*_i}$, and then map the rest of the vertices into $S\backslash \{t^*_1,\ldots, t^*_i\}$ bijectively.
    To show that this is indeed a homomorphism, notice first that $\{v_1,\ldots, v_{k-i-1},w_{t_0^*},\ldots, w_{t_i^*}\}$ is a partial edge in $F^{(\vec{t^*})}$.
    Therefore it remains to check that $\{w_{s_1},w_{s_2}\}$ and $\{v_{m},w_{s_1}\}$ are both in $\bigcup_{\vec{t}\in T}F^{(\vec{t})}$ for any $s_1\neq s_2\in S$ and $m\in[k-i-1]$.
    Indeed, if $s_1<s_2$ and $s_2 = t_j$ for some $\vec{t}\in T$, then $\{v_m, w_{s_1}, w_{s_2}\}$ is indeed a partial edge in $F^{(\vec{t})}$, which shows that both $\{w_{s_1},w_{s_2}\}$ and  $\{v_{m},w_{s_1}\}$ are partial edges in $F^{(\vec{t})}$ as well.
\end{proof}

We remark that \cref{thm:fix-k-big-r-full-edge} seems much stronger than \cref{cor:fix-k-big-r-full-edge}, though we do not see a clean way to extract a stronger statement from \cref{thm:fix-k-big-r-full-edge}.
We leave this as a potential future direction for interested readers.

With a completely different method, we can improve Mubayi's result in a slightly different way, and this is closer to what Sidorenko actually did in his paper \cite{Sid89} using hypergraph Lagrangian.
In that paper, Sidorenko showed that many extensions of partial $k$-graphs on $r+1$ vertices have Tur\'an density equal to $b(K_r^{(k)})$, as long as $r$ is at least some threshold $M_k$ that depends on $k$.
One special case related to our result is the $k$-graph $F_{r+1}^{(k,k-1)}$ that can be obtained as follows: consider the partial $k$-graph on $[r+1]$ spanned by the edges $\{[k-1]\cup i:i=k,\ldots, r+1\}$ and all the $2$-subsets of $[r+1]$, and then take the extension of the partial $k$-graph.
For example, $F_{k+1}^{(k,k-1)}$ is the tent $\Delta_{(k-1,1)}.$
Sidorenko's result is more general and relies on trees $T$ that satisfy the Erd\H{o}s--S\'os conjecture $\textup{ex}(T,n) \leq \frac{1}{2}(v(T)-2)n$, and we refer the readers to Sidorenko's original paper \cite{Sid89} for more details (also see \cite[Section 2]{Ste20} or \cite{TT22} for some families of trees where the Erd\H{o}s--S\'os conjecture is known to hold).

With a slightly different choice of partial forests, we can also prove that $\pi(F^{(k,k-1)}_{r+1}) = b(K^{(k)}_r)$ for sufficiently large $r$ with respect to $k$.
Our argument actually gives a more general statement: for any $s<k\leq r$, let $F^{(k,s)}_{r+1}$ be the extension of the partial $k$-graph spanned by $\{[s]\cup i:i=s+1,\ldots, r+1\}$ and all the $2$-subsets of $[r+1]$.
Then we obtain a sufficient condition for $\pi(F^{(k,s)}_{r+1}) = b(K_r^{(k)})$.

\begin{theorem}\label{thm:Sidorenko-type}
    Let $k,r,s$ be positive integers with $k\leq r$ and 
    \begin{equation}\label{eq:krs-relation}
        k-s\geq \sum_{i=1}^{s-1}\frac{i}{r-i}.
    \end{equation}
    Then $\pi(F^{(k,s)}_{r+1}) = b(K_r^{(k)})$.
\end{theorem}
\begin{proof}
    It is clear that $K_r^{(k)}$ is $F^{(k,s)}_{r+1}$-hom-free. Therefore, $\pi(F^{(k,s)}_{r+1})\geq b(K_r^{(k)})$ by \cref{thm:supersaturation}.
    
    To prove the other direction $\pi(F^{(k,s)}_{r+1})\leq b(K_r^{(k)})$,
    we may fix a $F^{(k,s)}_{r+1}$-hom-free $k$-graph $G$ and a random with uniform ordering $(X_1,\dots,X_k)$ on $G$.
    Let $x_1,\dots,x_k$ be the ratio sequence of $(X_1,\dots,X_k)$. We will solve for the maximum of $x_1\dots x_{k-1}$ under the constraints given by the following lemma.
    \begin{lemma}\label{lemma:sido-type-forest}
        For any integers $i,j$ with $i\in [k-s], i\leq j< k$, we have
        \[\frac{x_i}{r-k+i}\leq x_{j+1}-x_j.\]
    \end{lemma}
    \begin{proof}

    We will fix $i,j$ throughout this proof.
    As in what we did in \cref{sec:entropy}, we will temporarily fix an integer $N\in\NN$ that will later be taken to infinity.
    For any $t\in [N]$, we will define a partial forest $F^{(t)}$ on $V=\{v_1,\dots,v_{k-i},w_1,\dots,w_N\}$.
    The partial forest $F^{(t)}$ is spanned by the partial edges $\{v_1,\dots,v_{k-i},w_t\}$, $\{v_1,\dots,v_{k-j-1},w_m,w_t\}$ for every $m<t$, and $\{v_1,\dots,v_{k-j-1},w_m\}$ for every $m>t$.
    With the linear order $<$ given by $v_1<\dots<v_{k-i}<w_N<\cdots <w_1$, we know that $F^{(t)}$ is indeed a partial forest.
    We can compute the forest sequence with respect to the linear order as follows: each $v_m$ contributes one to $n_m$ for each $m\leq k-i$. For the contribution of $w_m$, if $m> t$ it contributes one to $n_{k-j}$; if $m=t$ it contributes one to $n_{k-i+1}$; otherwise it contributes one to $n_{k-j+1}$.
    Therefore the forest sequence $(n_1,\ldots, n_k)$ is $\vec{e}_1+\dots+\vec{e}_{k-i}+(N-t)\vec{e}_{k-j}+\vec{e}_{k-i+1}+(t-1)\vec{e}_{k-j+1}$, where $\vec{e}_1,\dots,\vec{e}_k$ are the vectors in the standard basis.
    Now let $(Y^{(t)}_v)_{v\in V}$ be the random homomorphism produced by \cref{lemma:sample-tree}.
    This gives
    \begin{align}\label{eq:sido-type-tree-eq1}
        \HH\left((Y^{(t)}_v)_{v\in V}\right) = (N+k-i)\HH(X_1)+\log_2\left(x_{i}\cdots x_{k}\cdot x_j^{t-1}x_{j+1}^{N-t}\right).
    \end{align}

    Now, we show that the random tuples $(Y^{(1)}_v)_{v\in V},\dots,(Y^{(N)}_v)_{v\in V}$ have $(r-k+i+1)$-wise disjoint supports.
    Note that, for any $t_1<\dots<t_{r-k+i+1}$, the extension of the union $\cup_{\ell=1}^{r-k+i+1} F^{(t_{\ell})}$ contains a homomorphic image of $F^{(k,k-i)}_{r+1}$, given by the partial edges $\{v_1,\dots,v_{k-i},w_{t_{\ell}}\}$ for $\ell\in [r-k+i+1]$ and $\{w_{t_{\ell'}},w_{t_{\ell}}\}$ for $1\leq \ell'<\ell\leq r-k+i+1$. Since $k-i\leq s$, this is also a homomorphic image of $F^{(k,s)}_{r+1}$. Thus, no sequence of vertices is in $\cap_{\ell=1}^{r-k+i+1}\supp ((Y^{(t_\ell)}_v)_{v\in V})$.
    
    Therefore we may now apply \cref{lemma:mix} with $a=r-k+i$.
    Suppose that $(Z_v)_{v\in V}$ is the resulting mixture of $(Y_v^{(t)})_{v\in V}$ for all $t\in [N]$.    
    Note that the partial edge $\{v_1,\dots,v_{k-i}\}$ is present in all partial forests, so by \cref{lemma:sample-tree} we know that $(Z_{v_1},\ldots, Z_{v_{k-i}})$ has the same distribution as $(X_{i+1},\ldots, X_k)$. Similarly, for each $m\in[N]$, since the partial edge $\{v_1,\dots,v_{k-j-1},w_m\}$ is present in all partial forests, we know that $(Z_{v_1},\ldots, Z_{v_{k-j-1}},Z_{w_m})$ has the same distribution as $(X_{j+1},\ldots, X_k)$.
    Hence
    \begin{align}
        \HH\left((Z_v)_{v\in V}\right)\leq& \HH(Z_{v_1},\ldots, Z_{v_{k-i}})+\sum_{m=1}^{N}\HH(Z_{w_m}\mid Z_{v_1},\ldots, Z_{v_{k-j-1}}) \nonumber\\
        =& \HH(X_{i+1},\ldots,X_k)+N\HH(X_{j+1}\mid X_{j+2},\ldots,X_k)\nonumber \\
        =& (N+k-i)\HH(X_1)+\log_2(x_{i+1}\cdots x_k\cdot x_{j+1}^N).\label{eq:sido-type-tree-eq2}
    \end{align}
    Thus \cref{lemma:mix,eq:sido-type-tree-eq1,eq:sido-type-tree-eq2} now give
    \[\sum_{t=1}^N x_{i}\cdots x_k\cdot  x_j^{t-1}x_{j+1}^{N-t}\cdot 2^{(N+k-i)\HH(X_1)}\leq (r-k+i)x_{i+1}\cdots x_k\cdot x_{j+1}^{N}\cdot 2^{(N+k-i)\HH(X_1)},\]
    and so
    \[\sum_{t=1}^Nx_ix_j^{t-1}x_{j+1}^{-t}\leq r-k+i.\]
    By rearranging and taking $N$ goes to infinity, we obtain
    \[\frac{x_i}{x_{j+1}}\cdot \frac{1}{1-\frac{x_j}{x_{j+1}}}=\sum_{t=1}^\infty \frac{x_i}{x_{j+1}}\left(\frac{x_j}{x_{j+1}}\right)^{t-1}\leq r-k+i,\]
    and the lemma follows.
    \end{proof}

Once again, to prove \cref{thm:Sidorenko-type}, we need to upper bound $x_1\cdots x_{k-1}$ given the inequalities in \cref{lemma:sido-type-forest}. We will show that the product is maximized when $x_i = (r-k+i)/r$ for each $i=1,\ldots, k-1$. We start with the following inequality similar to \cref{lemma:aux-ineq}.

\begin{lemma}\label{lemma:sido-type-aux-ineq}
    Suppose that $y_1,\ldots, y_t$ and $z$ are some non-negative real numbers.
    Then
    \[y_1\cdots y_t\leq \left(\sum_{i=1}^{t}\frac{y_i}{z+i}\right)^t\frac{(z+1)\cdots(z+t)}{t^t}.\]
\end{lemma}
\begin{proof}
        We will prove this by inducting on $t$. For $t=1$, the inequality is trivial.
        
        Assume the statement is true for $t-1$. From the inductive hypothesis and AM-GM inequality, we have
        \begin{align*}
            y_1\cdots y_t\leq& y_t\left(\sum_{i=1}^{t-1}\frac{y_i}{z+i}\right)^{t-1}\frac{(z+1)\cdots(z+t-1)}{(t-1)^{t-1}}\\
            = &\left(\frac{t-1}{z+t}y_t\right)\left(\sum_{i=1}^{t-1}\frac{y_i}{z+i}\right)^{t-1}\frac{(z+1)\cdots(z+t)}{(t-1)^{t}}\\
            \leq  &\left(\frac{t-1}{t}\sum_{i=1}^{t}\frac{y_i}{z+i}\right)^{t}\frac{(z+1)\cdots(z+t)}{(t-1)^{t}}\\
            =  &\left(\sum_{i=1}^{t}\frac{y_i}{z+i}\right)^t\frac{(z+1)\cdots(z+t)}{t^t}.\qedhere
        \end{align*}
\end{proof}

Now, by using this lemma with $t=k-1, y_i=x_i$ and $z=r-k$, it is sufficient to upper bound right hand side using the conditions from \cref{lemma:sido-type-forest}.
\begin{claim}\label{claim:sido-type-xk}
    We have
    \[\frac{x_1}{r-k+1}+\dots+\frac{x_{k-1}}{r-1}\leq \frac{k-1}{r}x_k.\]
\end{claim}
\begin{proof}
    Let $s'$ be the largest integer such that
    \[k-s'\geq \sum_{i=1}^{s'-1}\frac{i}{r-i}\]
    holds. In particular, we have $s\leq s'<k$.
    Set $c$ to be the real number such that
    \[\frac{k-1}{r}=(1-c)\frac{1}{r-s'}+\frac{1}{r-s'+1}+\dots+\frac{1}{r-1}.\]
    From the definition of $s'$, we have
    \[k-1\geq s'-1+\frac{s'-1}{r-s'+1}+\dots+\frac{1}{r-1}=\frac{r}{r-s'+1}+\dots+\frac{r}{r-1}\]
    and
    \[k-1< s'+\frac{s'}{r-s'}+\dots+\frac{1}{r-1}=\frac{r}{r-s'}+\dots+\frac{r}{r-1}.\]
    Therefore, $c\in (0,1]$. By replacing the coefficient of $x_k$ using the definition of $c$ and rearranging, we may rewrite the inequality we want to show as the following.
    \begin{align}
        &\frac{x_1}{r-k+1}+\dots+\frac{x_{k-s'-1}}{r-s'-1}+c\frac{x_{k-s'}}{r-s'}\nonumber\\
        \leq &(1-c)\frac{s'}{r-s'}\frac{x_k-x_{k-s'}}{s'}+\frac{s'-1}{r-s'+1}\frac{x_k-x_{k-s'+1}}{s'-1}+\dots+\frac{1}{r-1}\frac{x_k-x_{k-1}}{1}.\label{eq:sido-type-wts}
    \end{align}

    Note that \cref{lemma:sido-type-forest} implies that
    \[\frac{x_i}{r-k+i}\leq \frac{x_k-x_{k-j}}{j}\]
    holds for all $i\leq k-j$. Thus, to prove \cref{eq:sido-type-wts}, it is sufficient to check
    \[k-s'-1+c\leq (1-c)\frac{s'}{r-s'}+\frac{s'-1}{r-s'+1}+\dots+\frac{1}{r-1}.\]
    Actually, the equality holds because, by the choice of $c$, we have
    \begin{align*}
        k-s'-1+c=&(1-c)\frac{r}{r-s'}+\frac{r}{r-s'+1}+\frac{r}{r-s'+2}+\dots+\frac{r}{r-1}-s'+c\\
        =&(1-c)\frac{s'}{r-s'}+\frac{s'-1}{r-s'+1}+\frac{s'-2}{r-s'+2}+\dots+\frac{1}{r-1}.\qedhere
    \end{align*}
\end{proof}
By combining \cref{lemma:sido-type-aux-ineq,claim:sido-type-xk}, we get
\[x_1\dots x_{k-1}\leq \left(\frac{k-1}{r}x_k\right)^{k-1}\frac{(r-k+1)\cdots(r-1)}{(k-1)^{k-1}}=\frac{(r-k+1)\cdots(r-1)}{r^{k-1}}=b(K_r^{(k)}).\qedhere\]
\end{proof}

To give a sense of what the inequality in \cref{thm:Sidorenko-type} means, with some standard computation, we can show the following. 
If $r,k$ are growing positive integers such that $r = (C+o_{k\to\infty}(1))k$ for some $C\geq 1$, then the largest positive integer $s$ satisfying \cref{eq:krs-relation} is $(C(1-\exp(-C^{-1}))+o_{k\to\infty}(1))k$.
In a different regime where $s=k-d$ for some fixed positive integer $d$, we can get that the smallest positive integer $r$ satisfying the inequality is $((2d)^{-1}+o_{d;k\to\infty}(1))k^2$.
We include those computations in the appendix (\cref{prop:linear-r,prop:quadratic-r}).

We briefly remark that the threshold $M_k$ Sidorenko deduced on $r$ is the same as ours when $s=k-1$.
However, Sidorenko's argument works for a more general family of hypergraphs.
It is also possible that by modifying Sidorenko's argument appropriately, we may get a statement analogous to \cref{thm:Sidorenko-type} with the extra parameter $s$.

\section{Concluding remarks}\label{sec:conclusion}

\subsection{Exact result and stability}\label{subsec:stability}
In this paper, we mostly focus on the Tur\'an density rather than the Tur\'an number. 
However, we believe that with more work, it is possible to extract the exact Tur\'an number for sufficiently many vertices from our density Tur\'an theorems \cref{thm:main-tent,thm:Sidorenko-type} at least when we also forbid all homomorphic images.
More specifically, we believe that there are a corresponding stability results for \cref{thm:main-tent,thm:Sidorenko-type}, which is usually helpful to deduce the exact Tur\'an number for sufficiently many vertices.
Indeed, many exact results were deduced using stability results in a crucial way.
For some examples, we refer the readers to \cite{KS05,MP07, Pik08, Pik13, BIJ17, NY17, NY18, LMR23, San24}.

\subsection{Other extremizers}
All the Tur\'an results we are able to prove in this paper have blowups of $K^{(k)}_r$ as their asymptotic extremizers, and this is not a coincidence.
We find it much easier to construct partial forests that would give tight inequalities on the ratio sequences $x_1,\ldots, x_{k}$ with equality holding when $(X_1,\ldots,X_k)$ is a uniform oriented edge in $K^{(k)}_r$.
However, as mentioned in the introduction, many difficulties of hypergraph Tur\'an problems come from the potential complicated structures in the extremizers.
It would thus be more exciting if our method can be applied to problems with extremizers not as simple as $K^{(k)}_r$.

The first step would probably be to extend this to other Tur\'an problems where the extremizers are blowups of some other hypergraphs.
Two candidates are the complete bipartite $3$-graph $(A\sqcup B, E)$ where $E  = \binom{A}{2}\times B\cup A\times\binom{B}{2}$, and the complete oddly bipartite $k$-graph $(A\sqcup B,E)$ where $k$ is even, and $E$ is the $k$-edges $e$ such that $\abs{A\cap e}$ is odd.
Although they are not formally blowups of some smaller hypergraphs, one can think of the complete bipartite $3$-graphs as the blowups of $(\{1,2\},\{\{1,1,2\},\{1,2,2\}\})$, and the completely oddly bipartite $k$-graphs are the blowups of some $2$-vertex ``degenerate'' hypergraphs as well.

There are many known Tur\'an results where the two hypergraphs are (asymptotic) extremizers.
For example, a classical result of De Caen and F\"uredi \cite{DCF00} shows that the complete bipartite $3$-graph is an asymptotic extremizer for the Fano plane.
This was later extended by Mubayi--R\"odl \cite{MR02} and Baber--Talbot \cite{BT12}.
On the other hand, Keevash and Sudakov \cite{KS05} showed that the complete oddly bipartite $2k$-graph is the extremizer for the expanded triangle, the hypergraph with edges $\{1,\ldots, 2k\},\{k+1,\ldots, 3k\},\{1,\ldots,k,2k+1,\ldots,3k\}$.
A very recent breakthrough of Sankar \cite{San24} showed that the complete oddly bipartite $4$-graph is an asymptotic extremizer for tight cycles of sufficiently large length not divisible by $4$.

We are unable to construct any partial forests that give tight inequalities when $G$ is the complete bipartite $3$-graph.
For $G$ being complete oddly bipartite $k$-graphs, it is possible to construct such partial forests following the argument in \cref{thm:entropic-turan} and Sidorenko's \cite{Sid92} and Frankl's \cite{Fra90} ideas, which used auxiliary $2$-graphs to show that the Tur\'an densities of expanded triangles are $1/2$.
However, we have not found any other partial forests that use essentially different ideas.
It would be interesting to see if there are ways to obtain tight inequalities for those two candidates of $G$ in the hope that they would give rise to new Tur\'an results.

Let us close this discussion by mentioning that our method seems to capture a little structure in the conjectured extremizer for $K_4^{(3)-}$, the $3$-graph on $4$ vertices with $3$ edges.
Let $G_1$ be a $3$-graph on $6$ vertices with $10$ edges so that any $2$-subset is in exactly $2$ edges---it turns out that $G_1$ does exist and is unique up to isomorphism.
The \emph{iterated blowup} $G_m$ of $G_1$ is constructed inductively by replacing each vertex in $G_1$ with $G_{m-1}$.
Then $G_m$ is $K_4^{(3)-}$-free, and by taking $m$ to infinity, we get that $\pi(K_4^{(3)-})\geq \frac{2}{7}$.
This is a construction of Frankl and F\"uredi \cite{FF84}, and the construction is conjectured to be optimal.
The current best upper bound $\pi(K_4^{(3)-})\leq 0.2871$ is obtained by Baber and Talbot \cite{BT11} using flag algebra.
Though we cannot say anything new about the Tur\'an problem of $K_4^{(3)-}$ itself, our method seems to capture some structure in $G_1$.
Indeed, by the partial forests $F^{(i)}=([4],\{[3],[4]\backslash \{i\}\})$ for $i=1,2,3$, we can show that if $G$ is $K_4^{(3)-}$-free and $(X_1,X_2,X_3)$ is a random edge with uniform ordering on $G$, then
\[x_1 \eqdef 2^{\HH(X_1\mid X_2,X_3)-\HH(X_1)} \leq \frac{1}{3}.\]
This is indeed achieved when $(X_1,X_2,X_3)$ is a uniformly chosen oriented edge in $G_1$.

\subsection{Entropic spectral radius}
In \cref{sec:connection}, we showed that for any $k$-graph $G$, its spectral radius (i.e.\ the $k$-spectral radius) is related to the maximum of $\HH(X_2,\ldots, X_k\mid X_1)$ for symmetric distribution $(X_1,\ldots,X_k)$ on the oriented edges of $G$.
It would be interesting if this connection can be utilized to deduce some properties of spectral radius.
One possible candidate is a result of Kang, Liu and Shan \cite{KLS18} that showed that
\[\rho(G)\geq \left(\frac{1}{v(G)}\sum_{v\in V(G)}\deg(v)^{\frac{k}{k-1}}\right)^{\frac{k-1}{k}}\]
for any $k$-graph $G$, where $\rho(G)$ is the spectral radius of $G$.

\subsection{Entropic flag algebra}
As one may have observed, many upper bounds on Tur\'an densities, especially for those that are still open, were obtained using flag algebra.
Such upper bounds using flag algebra, roughly speaking, are obtained via carefully chosen sum-of-squares inequalities, enumeration of possible small configurations, and numerical computation of positive semidefinite programs.
See \cite{Raz13} for a more detailed discussion of the method.

The inequalities obtained using our argument seem to be really different from the inequalities obtained by sum-of-squares.
This suggests a possibility that maybe the flag algebra bounds can be improved with this new idea and some enumeration of possible partial forests to use in the argument.
However, aside from the time complexity enumerating through the possible partial forests, there seem to be several technicalities to overcome for this to work.
The first is that in most of our proofs, we need to look at infinitely many partial forests in order to get a tight bound.
In addition, the inequalities we get, unlike the ones in flag-algebraic arguments, are highly non-linear.
However, if we are just aiming for some numerical upper bound that is close to the truth, then hopefully finite but sufficiently many partial forests together with an approximation of the supremum of $x_1\cdots x_{k-1}$ subject to the inequalities would be enough.

The most serious issue is probably that there has not been a framework for automated entropic computation.
So far, the flag-algebraic tools are developed to keep track of the homomorphism densities of labeled graphs.
Unfortunately, it seems that all our arguments for hypergraph Tur\'an problems cannot be rephrased using homomorphism densities as we also crucially use the marginal distributions of the random homomorphisms sampled by \cref{lemma:sample-tree}.
It would thus be necessary to come up with an ``entropic flag algebra'' framework and implement corresponding software to execute the idea in this subsection.
We refer the readers to \cite{CY24-2} for another entropic argument that motivates this idea of ``entropic flag algebra''.
\section*{Acknowledgement}
The project was motivated when the first author was visiting Hong Liu at Institute for Basic Science, and the first author would like to thank his hospitality.
We would also like to thank Ryan Alweiss and Freddie Manners for discussions during the early stage of this project, Dhruv Mubayi and Maya Sankar for pointing us to references for hypergraph Tur\'an problems, Yongtao Li for pointing us to references for spectral Tur\'an problems, and Noga Alon for pointing us to other useful references.
In addition, we would like to thank Zeev Dvir, Xiaoyu He, Cosmin Pohoata and Maya Sankar for helpful comments on an earlier draft.
Last but not least, we thank the anonymous reviewer for detailed and helpful comments.

The second author is supported by NSF grant DMS 2246682.

\bibliographystyle{amsplain0}
\bibliography{ref_joints}

@article {FK98,
    AUTHOR = {Friedgut, Ehud and Kahn, Jeff},
     TITLE = {On the number of copies of one hypergraph in another},
   JOURNAL = {Israel J. Math.},
  FJOURNAL = {Israel Journal of Mathematics},
    VOLUME = {105},
      YEAR = {1998},
     PAGES = {251--256},
      ISSN = {0021-2172},
   MRCLASS = {05C65},
  MRNUMBER = {1639767},
MRREVIEWER = {Nigel Martin},
       DOI = {10.1007/BF02780332},
       URL = {https://doi.org/10.1007/BF02780332},
}

@article {CGFS86,
    AUTHOR = {Chung, F. R. K. and Graham, R. L. and Frankl, P. and Shearer,
              J. B.},
     TITLE = {Some intersection theorems for ordered sets and graphs},
   JOURNAL = {J. Combin. Theory Ser. A},
  FJOURNAL = {Journal of Combinatorial Theory. Series A},
    VOLUME = {43},
      YEAR = {1986},
    NUMBER = {1},
     PAGES = {23--37},
      ISSN = {0097-3165},
   MRCLASS = {05A05 (05C65)},
  MRNUMBER = {859293},
MRREVIEWER = {Zolt\'{a}n F\"{u}redi},
       DOI = {10.1016/0097-3165(86)90019-1},
       URL = {https://doi.org/10.1016/0097-3165(86)90019-1},
}

@book {AS00,
    AUTHOR = {Alon, Noga and Spencer, Joel H.},
     TITLE = {The probabilistic method},
    SERIES = {Wiley-Interscience Series in Discrete Mathematics and
              Optimization},
   EDITION = {Second},
      NOTE = {With an appendix on the life and work of Paul Erd\H{o}s},
 PUBLISHER = {Wiley-Interscience [John Wiley \& Sons], New York},
      YEAR = {2000},
     PAGES = {xviii+301},
      ISBN = {0-471-37046-0},
   MRCLASS = {60-02 (05C80 60C05 60F99 60G42)},
  MRNUMBER = {1885388},
MRREVIEWER = {Bert Fristedt},
       DOI = {10.1002/0471722154},
       URL = {https://doi.org/10.1002/0471722154},
}

@article {CY24,
    AUTHOR = {Chao, Ting-Wei and Yu, Hung-Hsun Hans},
     TITLE = {Kruskal--{K}atona-type problems via the entropy method},
   JOURNAL = {J. Combin. Theory Ser. B},
  FJOURNAL = {Journal of Combinatorial Theory. Series B},
    VOLUME = {169},
      YEAR = {2024},
     PAGES = {480--506},
      ISSN = {0095-8956,1096-0902},
   MRCLASS = {99-06},
  MRNUMBER = {4789260},
       DOI = {10.1016/j.jctb.2024.08.003},
       URL = {https://doi.org/10.1016/j.jctb.2024.08.003},
}

@misc{CY24-2,
      title={A Purely Entropic Approach to the Rainbow Triangle Problem}, 
      author={Ting-Wei Chao and Hung-Hsun Hans Yu},
      Note = {\arXiv{2407.14084}}
}

@phdthesis{Fitch18,
  title={Applications of entropy to extremal problems},
  author={Fitch, Matthew},
  year={2018},
  school={University of Warwick}
}

@incollection {Kru63,
    AUTHOR = {Kruskal, Joseph B.},
     TITLE = {The number of simplices in a complex},
 BOOKTITLE = {Mathematical optimization techniques},
     PAGES = {251--278},
 PUBLISHER = {Univ. California Press, Berkeley-Los Angeles, Calif.},
      YEAR = {1963},
   MRCLASS = {05.10},
  MRNUMBER = {154827},
MRREVIEWER = {John\ Riordan},
}

@incollection {Kat68,
    AUTHOR = {Katona, G.},
     TITLE = {A theorem of finite sets},
 BOOKTITLE = {Theory of {G}raphs ({P}roc. {C}olloq., {T}ihany, 1966)},
     PAGES = {187--207},
 PUBLISHER = {Academic Press, New York-London},
      YEAR = {1968},
   MRCLASS = {05A05},
  MRNUMBER = {290982},
MRREVIEWER = {P.\ Erd\H{o}s},
}

@article {Turan41,
    AUTHOR = {Tur\'an, Paul},
     TITLE = {Eine {E}xtremalaufgabe aus der {G}raphentheorie},
   JOURNAL = {Mat. Fiz. Lapok},
  FJOURNAL = {Matematikai \'es Fizikai Lapok},
    VOLUME = {48},
      YEAR = {1941},
     PAGES = {436--452},
      ISSN = {0302-7317},
   MRCLASS = {56.0X},
  MRNUMBER = {18405},
MRREVIEWER = {P.\ Erd\H os},
}

@article {ES46,
    AUTHOR = {Erd\"os, P. and Stone, A. H.},
     TITLE = {On the structure of linear graphs},
   JOURNAL = {Bull. Amer. Math. Soc.},
  FJOURNAL = {Bulletin of the American Mathematical Society},
    VOLUME = {52},
      YEAR = {1946},
     PAGES = {1087--1091},
      ISSN = {0002-9904},
   MRCLASS = {56.0X},
  MRNUMBER = {18807},
MRREVIEWER = {H.\ S. M. Coxeter},
       DOI = {10.1090/S0002-9904-1946-08715-7},
       URL = {https://doi.org/10.1090/S0002-9904-1946-08715-7},
}

@article {Raz10,
    AUTHOR = {Razborov, Alexander A.},
     TITLE = {On 3-hypergraphs with forbidden 4-vertex configurations},
   JOURNAL = {SIAM J. Discrete Math.},
  FJOURNAL = {SIAM Journal on Discrete Mathematics},
    VOLUME = {24},
      YEAR = {2010},
    NUMBER = {3},
     PAGES = {946--963},
      ISSN = {0895-4801,1095-7146},
   MRCLASS = {05C65 (90C22)},
  MRNUMBER = {2680226},
MRREVIEWER = {Peter\ D.\ Johnson, Jr.},
       DOI = {10.1137/090747476},
       URL = {https://doi.org/10.1137/090747476},
}

@article {F-RV13,
    AUTHOR = {Falgas-Ravry, Victor and Vaughan, Emil R.},
     TITLE = {Applications of the semi-definite method to the {T}ur\'an
              density problem for 3-graphs},
   JOURNAL = {Combin. Probab. Comput.},
  FJOURNAL = {Combinatorics, Probability and Computing},
    VOLUME = {22},
      YEAR = {2013},
    NUMBER = {1},
     PAGES = {21--54},
      ISSN = {0963-5483,1469-2163},
   MRCLASS = {05D05 (05C65)},
  MRNUMBER = {3002572},
       DOI = {10.1017/S0963548312000508},
       URL = {https://doi.org/10.1017/S0963548312000508},
}

@article {Kos82,
    AUTHOR = {Kostochka, A. V.},
     TITLE = {A class of constructions for {T}ur\'an's {$(3,\,4)$}-problem},
   JOURNAL = {Combinatorica},
  FJOURNAL = {Combinatorica. An International Journal of the J\'anos Bolyai
              Mathematical Society},
    VOLUME = {2},
      YEAR = {1982},
    NUMBER = {2},
     PAGES = {187--192},
      ISSN = {0209-9683},
   MRCLASS = {05C35 (05C65)},
  MRNUMBER = {685045},
MRREVIEWER = {K.\ Vesztergombi},
       DOI = {10.1007/BF02579317},
       URL = {https://doi.org/10.1007/BF02579317},
}

@article {F-D-F88,
    AUTHOR = {Fon-Der-Flaass, D. G.},
     TITLE = {A method for constructing {$(3,4)$}-graphs},
   JOURNAL = {Mat. Zametki},
  FJOURNAL = {Akademiya Nauk SSSR. Matematicheskie Zametki},
    VOLUME = {44},
      YEAR = {1988},
    NUMBER = {4},
     PAGES = {546--550, 559},
      ISSN = {0025-567X},
   MRCLASS = {05C35},
  MRNUMBER = {975195},
MRREVIEWER = {Pavel\ Tomasta},
       DOI = {10.1007/BF01158925},
       URL = {https://doi.org/10.1007/BF01158925},
}

@article {Erdos70,
    AUTHOR = {Erd\H{o}s, P\'al},
     TITLE = {On the graph theorem of {T}ur\'an},
   JOURNAL = {Mat. Lapok},
  FJOURNAL = {Matematikai Lapok. Bolyai J\'anos Matematikai T\'arsulat},
    VOLUME = {21},
      YEAR = {1970},
     PAGES = {249--251},
      ISSN = {0025-519X},
   MRCLASS = {05C35},
  MRNUMBER = {307975},
MRREVIEWER = {\v S.\ Zn\'am},
}

@book {AZ18,
    AUTHOR = {Aigner, Martin and Ziegler, G\"unter M.},
     TITLE = {Proofs from {T}he {B}ook},
   EDITION = {Sixth},
 PUBLISHER = {Springer, Berlin},
      YEAR = {2018},
     PAGES = {viii+326},
      ISBN = {978-3-662-57264-1; 978-3-662-57265-8},
   MRCLASS = {00A05},
  MRNUMBER = {3823190},
       DOI = {10.1007/978-3-662-57265-8},
       URL = {https://doi.org/10.1007/978-3-662-57265-8},
}

@techreport{Caro79,
  title={New results on the independence number},
  author={Caro, Yair},
  year={1979},
  institution={Technical Report, Tel-Aviv University}
}

@article {MS65,
    AUTHOR = {Motzkin, T. S. and Straus, E. G.},
     TITLE = {Maxima for graphs and a new proof of a theorem of {T}ur\'an},
   JOURNAL = {Canadian J. Math.},
  FJOURNAL = {Canadian Journal of Mathematics. Journal Canadien de
              Math\'ematiques},
    VOLUME = {17},
      YEAR = {1965},
     PAGES = {533--540},
      ISSN = {0008-414X,1496-4279},
   MRCLASS = {05.50},
  MRNUMBER = {175813},
MRREVIEWER = {F.\ Harary},
       DOI = {10.4153/CJM-1965-053-6},
       URL = {https://doi.org/10.4153/CJM-1965-053-6},
}

@article {Sha48,
    AUTHOR = {Shannon, C. E.},
     TITLE = {A mathematical theory of communication},
   JOURNAL = {Bell System Tech. J.},
  FJOURNAL = {The Bell System Technical Journal},
    VOLUME = {27},
      YEAR = {1948},
     PAGES = {379--423, 623--656},
      ISSN = {0005-8580},
   MRCLASS = {60.0X},
  MRNUMBER = {26286},
MRREVIEWER = {J.\ L.\ Doob},
       DOI = {10.1002/j.1538-7305.1948.tb01338.x},
       URL = {https://doi.org/10.1002/j.1538-7305.1948.tb01338.x},
}

@article {Rad97,
    AUTHOR = {Radhakrishnan, Jaikumar},
     TITLE = {An entropy proof of {B}regman's theorem},
   JOURNAL = {J. Combin. Theory Ser. A},
  FJOURNAL = {Journal of Combinatorial Theory. Series A},
    VOLUME = {77},
      YEAR = {1997},
    NUMBER = {1},
     PAGES = {161--164},
      ISSN = {0097-3165,1096-0899},
   MRCLASS = {15A15},
  MRNUMBER = {1426744},
       DOI = {10.1006/jcta.1996.2727},
       URL = {https://doi.org/10.1006/jcta.1996.2727},
}

@misc{Gil22,
      title={A constant lower bound for the union-closed sets conjecture}, 
      author={Justin Gilmer},
      Note = {\arXiv{2211.09055}}
}

@misc{GGMT24,
      title={On a conjecture of Marton}, 
      author={W. T. Gowers and Ben Green and Freddie Manners and Terence Tao},
      Note = {to appear on Ann. of Math.}
}

@misc{Sze15,
      title={An information theoretic approach to {S}idorenko's conjecture}, 
      author={Bal{\'a}zs Szegedy},
      Note = {\arXiv{1406.6738}}
}

@misc{LS11,
      title={On the logarithimic calculus and {S}idorenko's conjecture}, 
      author={J.L. Xiang Li and Bal{\'a}zs Szegedy},
      Note = {\arXiv{1107.1153}}
}

@misc{CKLL18-2,
      title={Sidorenko's conjecture for higher tree decompositions}, 
      author={David Conlon and Jeong Han Kim and Choongbum Lee and Joonkyung Lee},
      Note = {\arXiv{1805.02238}}
}

@misc{San24,
      title={The {Tur\'{a}n} Density of 4-Uniform Tight Cycles}, 
      author={Maya Sankar},
      Note = {\arXiv{2411.01782}}
}

@article {CKLL18-1,
    AUTHOR = {Conlon, David and Kim, Jeong Han and Lee, Choongbum and Lee,
              Joonkyung},
     TITLE = {Some advances on {S}idorenko's conjecture},
   JOURNAL = {J. Lond. Math. Soc. (2)},
  FJOURNAL = {Journal of the London Mathematical Society. Second Series},
    VOLUME = {98},
      YEAR = {2018},
    NUMBER = {3},
     PAGES = {593--608},
      ISSN = {0024-6107,1469-7750},
   MRCLASS = {05C35 (05C60 05D40)},
  MRNUMBER = {3893193},
MRREVIEWER = {J\'ozsef\ Balogh},
       DOI = {10.1112/jlms.12142},
       URL = {https://doi.org/10.1112/jlms.12142},
}

@phdthesis{Par14,
  title={On {S}idorenko’s conjecture},
  author={Parczyk, Olaf},
  year={2014},
  school={Master’s thesis, Freie Universit{\"a}t, Berlin}
}

@article {CL17,
    AUTHOR = {Conlon, David and Lee, Joonkyung},
     TITLE = {Finite reflection groups and graph norms},
   JOURNAL = {Adv. Math.},
  FJOURNAL = {Advances in Mathematics},
    VOLUME = {315},
      YEAR = {2017},
     PAGES = {130--165},
      ISSN = {0001-8708,1090-2082},
   MRCLASS = {05C60 (05C50 05C70 05C80 20D99)},
  MRNUMBER = {3667583},
MRREVIEWER = {Juanjo\ Ru\'e},
       DOI = {10.1016/j.aim.2017.05.009},
       URL = {https://doi.org/10.1016/j.aim.2017.05.009},
}

@article {GLLV22,
    AUTHOR = {Grzesik, Andrzej and Lee, Joonkyung and Lidick\'y, Bernard and
              Volec, Jan},
     TITLE = {On tripartite common graphs},
   JOURNAL = {Combin. Probab. Comput.},
  FJOURNAL = {Combinatorics, Probability and Computing},
    VOLUME = {31},
      YEAR = {2022},
    NUMBER = {5},
     PAGES = {907--923},
      ISSN = {0963-5483,1469-2163},
   MRCLASS = {05C35 (05C55)},
  MRNUMBER = {4472294},
MRREVIEWER = {Kiyoshi\ Yoshimoto},
       DOI = {10.1017/s0963548322000074},
       URL = {https://doi.org/10.1017/s0963548322000074},
}

@article {BMN24,
    AUTHOR = {Behague, Natalie and Morrison, Natasha and Noel, Jonathan A.},
     TITLE = {Off-diagonal commonality of graphs via entropy},
   JOURNAL = {SIAM J. Discrete Math.},
  FJOURNAL = {SIAM Journal on Discrete Mathematics},
    VOLUME = {38},
      YEAR = {2024},
    NUMBER = {3},
     PAGES = {2335--2360},
      ISSN = {0895-4801,1095-7146},
   MRCLASS = {05C55 (05C35)},
  MRNUMBER = {4790926},
       DOI = {10.1137/23M1625342},
       URL = {https://doi.org/10.1137/23M1625342},
}

@article {Lee21,
    AUTHOR = {Lee, Joonkyung},
     TITLE = {On some graph densities in locally dense graphs},
   JOURNAL = {Random Structures Algorithms},
  FJOURNAL = {Random Structures \& Algorithms},
    VOLUME = {58},
      YEAR = {2021},
    NUMBER = {2},
     PAGES = {322--344},
      ISSN = {1042-9832,1098-2418},
   MRCLASS = {05C60 (05C35 05C55)},
  MRNUMBER = {4201799},
       DOI = {10.1002/rsa.20974},
       URL = {https://doi.org/10.1002/rsa.20974},
}

@article{Kee11,
  title={Hypergraph {Tur\'an} problems},
  author={Keevash, Peter},
  journal={Surveys in combinatorics},
  volume={392},
  pages={83--140},
  year={2011},
  publisher={Cambridge University Press Cambridge}
}

@article {FR84,
    AUTHOR = {Frankl, P. and R\"odl, V.},
     TITLE = {Hypergraphs do not jump},
   JOURNAL = {Combinatorica},
  FJOURNAL = {Combinatorica. An International Journal of the J\'anos Bolyai
              Mathematical Society},
    VOLUME = {4},
      YEAR = {1984},
    NUMBER = {2-3},
     PAGES = {149--159},
      ISSN = {0209-9683},
   MRCLASS = {05C65},
  MRNUMBER = {771722},
MRREVIEWER = {Noga\ Alon},
       DOI = {10.1007/BF02579215},
       URL = {https://doi.org/10.1007/BF02579215},
}

@article {BT11,
    AUTHOR = {Baber, Rahil and Talbot, John},
     TITLE = {Hypergraphs do jump},
   JOURNAL = {Combin. Probab. Comput.},
  FJOURNAL = {Combinatorics, Probability and Computing},
    VOLUME = {20},
      YEAR = {2011},
    NUMBER = {2},
     PAGES = {161--171},
      ISSN = {0963-5483,1469-2163},
   MRCLASS = {05C65 (05C35)},
  MRNUMBER = {2769186},
MRREVIEWER = {Yi\ Zhao},
       DOI = {10.1017/S0963548310000222},
       URL = {https://doi.org/10.1017/S0963548310000222},
}

@article {Mub06,
    AUTHOR = {Mubayi, Dhruv},
     TITLE = {A hypergraph extension of {T}ur\'an's theorem},
   JOURNAL = {J. Combin. Theory Ser. B},
  FJOURNAL = {Journal of Combinatorial Theory. Series B},
    VOLUME = {96},
      YEAR = {2006},
    NUMBER = {1},
     PAGES = {122--134},
      ISSN = {0095-8956,1096-0902},
   MRCLASS = {05C65 (05C35)},
  MRNUMBER = {2185983},
MRREVIEWER = {Peter\ Keevash},
       DOI = {10.1016/j.jctb.2005.06.013},
       URL = {https://doi.org/10.1016/j.jctb.2005.06.013},
}

@article {Sid89,
    AUTHOR = {Sidorenko, A. F.},
     TITLE = {Asymptotic solution for a new class of forbidden {$r$}-graphs},
   JOURNAL = {Combinatorica},
  FJOURNAL = {Combinatorica. An International Journal on Combinatorics and
              the Theory of Computing},
    VOLUME = {9},
      YEAR = {1989},
    NUMBER = {2},
     PAGES = {207--215},
      ISSN = {0209-9683},
   MRCLASS = {05C35 (05C65)},
  MRNUMBER = {1030374},
MRREVIEWER = {Ioan\ Tomescu},
       DOI = {10.1007/BF02124681},
       URL = {https://doi.org/10.1007/BF02124681},
}

@article {Bol74,
    AUTHOR = {Bollob\'as, B\'ela},
     TITLE = {Three-graphs without two triples whose symmetric difference is
              contained in a third},
   JOURNAL = {Discrete Math.},
  FJOURNAL = {Discrete Mathematics},
    VOLUME = {8},
      YEAR = {1974},
     PAGES = {21--24},
      ISSN = {0012-365X,1872-681X},
   MRCLASS = {05C99},
  MRNUMBER = {345869},
       DOI = {10.1016/0012-365X(74)90105-8},
       URL = {https://doi.org/10.1016/0012-365X(74)90105-8},
}

@article {FF83,
    AUTHOR = {Frankl, Peter and F\"uredi, Zolt\'an},
     TITLE = {A new generalization of the {E}rd{\H{o}}s-{K}o-{R}ado theorem},
   JOURNAL = {Combinatorica},
  FJOURNAL = {Combinatorica. An International Journal of the J\'anos Bolyai
              Mathematical Society},
    VOLUME = {3},
      YEAR = {1983},
    NUMBER = {3-4},
     PAGES = {341--349},
      ISSN = {0209-9683},
   MRCLASS = {05A05 (05B30)},
  MRNUMBER = {729787},
MRREVIEWER = {\c Serban\ Buze\c teanu},
       DOI = {10.1007/BF02579190},
       URL = {https://doi.org/10.1007/BF02579190},
}

@article {BT12,
    AUTHOR = {Baber, Rahil and Talbot, John},
     TITLE = {New {T}ur\'an densities for 3-graphs},
   JOURNAL = {Electron. J. Combin.},
  FJOURNAL = {Electronic Journal of Combinatorics},
    VOLUME = {19},
      YEAR = {2012},
    NUMBER = {2},
     PAGES = {Paper 22, 21},
      ISSN = {1077-8926},
   MRCLASS = {05D05 (05C65)},
  MRNUMBER = {2928637},
MRREVIEWER = {Yi\ Zhao},
       DOI = {10.37236/2360},
       URL = {https://doi.org/10.37236/2360},
}

@article {Pik08,
    AUTHOR = {Pikhurko, Oleg},
     TITLE = {An exact {T}ur\'an result for the generalized triangle},
   JOURNAL = {Combinatorica},
  FJOURNAL = {Combinatorica. An International Journal on Combinatorics and
              the Theory of Computing},
    VOLUME = {28},
      YEAR = {2008},
    NUMBER = {2},
     PAGES = {187--208},
      ISSN = {0209-9683,1439-6912},
   MRCLASS = {05C65 (05C35 05D05)},
  MRNUMBER = {2399018},
MRREVIEWER = {J\'ozsef\ Balogh},
       DOI = {10.1007/s00493-008-2187-2},
       URL = {https://doi.org/10.1007/s00493-008-2187-2},
}

@article{Wei81,
  title={A lower bound on the stability number of a simple graph},
  author={Wei, Victor K},
  year={1981},
  publisher={Bell Laboratories Technical Memorandum Murray Hill, NJ, USA},
  JOURNAL = {Bell Lab. Tech. Memor.},
  FJOURNAL = {Bell Laboratories Technical Memorandum Murray Hill, NJ, USA},
}

@article {FF89,
    AUTHOR = {Frankl, P. and F\"uredi, Z.},
     TITLE = {Extremal problems whose solutions are the blowups of the small
              {W}itt-designs},
   JOURNAL = {J. Combin. Theory Ser. A},
  FJOURNAL = {Journal of Combinatorial Theory. Series A},
    VOLUME = {52},
      YEAR = {1989},
    NUMBER = {1},
     PAGES = {129--147},
      ISSN = {0097-3165,1096-0899},
   MRCLASS = {05A05 (05C35 05C65)},
  MRNUMBER = {1008165},
MRREVIEWER = {J.\ R.\ Griggs},
       DOI = {10.1016/0097-3165(89)90067-8},
       URL = {https://doi.org/10.1016/0097-3165(89)90067-8},
}

@article {Qi13,
    AUTHOR = {Qi, Liqun},
     TITLE = {Symmetric nonnegative tensors and copositive tensors},
   JOURNAL = {Linear Algebra Appl.},
  FJOURNAL = {Linear Algebra and its Applications},
    VOLUME = {439},
      YEAR = {2013},
    NUMBER = {1},
     PAGES = {228--238},
      ISSN = {0024-3795,1873-1856},
   MRCLASS = {15A69 (15A18)},
  MRNUMBER = {3045233},
MRREVIEWER = {Minru\ Bai},
       DOI = {10.1016/j.laa.2013.03.015},
       URL = {https://doi.org/10.1016/j.laa.2013.03.015},
}

@article {Nik02,
    AUTHOR = {Nikiforov, V.},
     TITLE = {Some inequalities for the largest eigenvalue of a graph},
   JOURNAL = {Combin. Probab. Comput.},
  FJOURNAL = {Combinatorics, Probability and Computing},
    VOLUME = {11},
      YEAR = {2002},
    NUMBER = {2},
     PAGES = {179--189},
      ISSN = {0963-5483,1469-2163},
   MRCLASS = {05C50},
  MRNUMBER = {1888908},
       DOI = {10.1017/S0963548301004928},
       URL = {https://doi.org/10.1017/S0963548301004928},
}

@article {Wil86,
    AUTHOR = {Wilf, Herbert S.},
     TITLE = {Spectral bounds for the clique and independence numbers of
              graphs},
   JOURNAL = {J. Combin. Theory Ser. B},
  FJOURNAL = {Journal of Combinatorial Theory. Series B},
    VOLUME = {40},
      YEAR = {1986},
    NUMBER = {1},
     PAGES = {113--117},
      ISSN = {0095-8956,1096-0902},
   MRCLASS = {05C50 (05C99)},
  MRNUMBER = {830598},
MRREVIEWER = {Qiao\ Li},
       DOI = {10.1016/0095-8956(86)90069-9},
       URL = {https://doi.org/10.1016/0095-8956(86)90069-9},
}

@article {NY17,
    AUTHOR = {Norin, S. and Yepremyan, L.},
     TITLE = {{Tur\'an} number of generalized triangles},
   JOURNAL = {J. Combin. Theory Ser. A},
  FJOURNAL = {Journal of Combinatorial Theory. Series A},
    VOLUME = {146},
      YEAR = {2017},
     PAGES = {312--343},
      ISSN = {0097-3165,1096-0899},
   MRCLASS = {05C35 (05C65)},
  MRNUMBER = {3574234},
MRREVIEWER = {Yi\ Zhao},
       DOI = {10.1016/j.jcta.2016.09.003},
       URL = {https://doi.org/10.1016/j.jcta.2016.09.003},
}

@article {Pik13,
    AUTHOR = {Pikhurko, Oleg},
     TITLE = {Exact computation of the hypergraph {T}ur\'an function for
              expanded complete 2-graphs},
   JOURNAL = {J. Combin. Theory Ser. B},
  FJOURNAL = {Journal of Combinatorial Theory. Series B},
    VOLUME = {103},
      YEAR = {2013},
    NUMBER = {2},
     PAGES = {220--225},
      ISSN = {0095-8956,1096-0902},
   MRCLASS = {05C65 (05C35)},
  MRNUMBER = {3018066},
MRREVIEWER = {Guantao\ Chen},
       DOI = {10.1016/j.jctb.2012.09.005},
       URL = {https://doi.org/10.1016/j.jctb.2012.09.005},
}

@article {BIJ17,
    AUTHOR = {Brandt, Axel and Irwin, David and Jiang, Tao},
     TITLE = {Stability and {T}ur\'an numbers of a class of hypergraphs via
              {L}agrangians},
   JOURNAL = {Combin. Probab. Comput.},
  FJOURNAL = {Combinatorics, Probability and Computing},
    VOLUME = {26},
      YEAR = {2017},
    NUMBER = {3},
     PAGES = {367--405},
      ISSN = {0963-5483,1469-2163},
   MRCLASS = {05C65 (05C35)},
  MRNUMBER = {3628909},
MRREVIEWER = {Ioan\ Tomescu},
       DOI = {10.1017/S0963548316000444},
       URL = {https://doi.org/10.1017/S0963548316000444},
}

@article {NY18,
    AUTHOR = {Norin, Sergey and Yepremyan, Liana},
     TITLE = {{Tur\'an} numbers of extensions},
   JOURNAL = {J. Combin. Theory Ser. A},
  FJOURNAL = {Journal of Combinatorial Theory. Series A},
    VOLUME = {155},
      YEAR = {2018},
     PAGES = {476--492},
      ISSN = {0097-3165,1096-0899},
   MRCLASS = {05C35 (05C65)},
  MRNUMBER = {3741438},
MRREVIEWER = {Yi\ Zhao},
       DOI = {10.1016/j.jcta.2017.08.004},
       URL = {https://doi.org/10.1016/j.jcta.2017.08.004},
}

@article {LMR23,
    AUTHOR = {Liu, Xizhi and Mubayi, Dhruv and Reiher, Christian},
     TITLE = {A unified approach to hypergraph stability},
   JOURNAL = {J. Combin. Theory Ser. B},
  FJOURNAL = {Journal of Combinatorial Theory. Series B},
    VOLUME = {158},
      YEAR = {2023},
     PAGES = {36--62},
      ISSN = {0095-8956,1096-0902},
   MRCLASS = {05C35 (05C65)},
  MRNUMBER = {4484827},
MRREVIEWER = {Jian\ Wang},
       DOI = {10.1016/j.jctb.2022.08.008},
       URL = {https://doi.org/10.1016/j.jctb.2022.08.008},
}

@article {DCF00,
    AUTHOR = {De Caen, Dominique and F\"uredi, Zolt\'an},
     TITLE = {The maximum size of 3-uniform hypergraphs not containing a
              {F}ano plane},
   JOURNAL = {J. Combin. Theory Ser. B},
  FJOURNAL = {Journal of Combinatorial Theory. Series B},
    VOLUME = {78},
      YEAR = {2000},
    NUMBER = {2},
     PAGES = {274--276},
      ISSN = {0095-8956,1096-0902},
   MRCLASS = {05C65 (05B25)},
  MRNUMBER = {1750899},
       DOI = {10.1006/jctb.1999.1938},
       URL = {https://doi.org/10.1006/jctb.1999.1938},
}

@article {MR02,
    AUTHOR = {Mubayi, Dhruv and R\"odl, Vojt\^ech},
     TITLE = {On the {T}ur\'an number of triple systems},
   JOURNAL = {J. Combin. Theory Ser. A},
  FJOURNAL = {Journal of Combinatorial Theory. Series A},
    VOLUME = {100},
      YEAR = {2002},
    NUMBER = {1},
     PAGES = {136--152},
      ISSN = {0097-3165,1096-0899},
   MRCLASS = {05C35},
  MRNUMBER = {1932073},
MRREVIEWER = {Sergei\ L.\ Bezrukov},
       DOI = {10.1006/jcta.2002.3284},
       URL = {https://doi.org/10.1006/jcta.2002.3284},
}

@article {KS05,
    AUTHOR = {Keevash, Peter and Sudakov, Benny},
     TITLE = {On a hypergraph {T}ur\'an problem of {F}rankl},
   JOURNAL = {Combinatorica},
  FJOURNAL = {Combinatorica. An International Journal on Combinatorics and
              the Theory of Computing},
    VOLUME = {25},
      YEAR = {2005},
    NUMBER = {6},
     PAGES = {673--706},
      ISSN = {0209-9683,1439-6912},
   MRCLASS = {05C65 (05C35 05D05 05E35)},
  MRNUMBER = {2199431},
MRREVIEWER = {Yi\ Zhao},
       DOI = {10.1007/s00493-005-0042-2},
       URL = {https://doi.org/10.1007/s00493-005-0042-2},
}

@article {Sid92,
    AUTHOR = {Sidorenko, Alexander},
     TITLE = {Asymptotic solution of the {T}ur\'an problem for some
              hypergraphs},
   JOURNAL = {Graphs Combin.},
  FJOURNAL = {Graphs and Combinatorics},
    VOLUME = {8},
      YEAR = {1992},
    NUMBER = {2},
     PAGES = {199--201},
      ISSN = {0911-0119,1435-5914},
   MRCLASS = {99-04},
  MRNUMBER = {1554360},
       DOI = {10.1007/BF02350635},
       URL = {https://doi.org/10.1007/BF02350635},
}

@article {Fra90,
    AUTHOR = {Frankl, P.},
     TITLE = {Asymptotic solution of a {T}ur\'an-type problem},
   JOURNAL = {Graphs Combin.},
  FJOURNAL = {Graphs and Combinatorics},
    VOLUME = {6},
      YEAR = {1990},
    NUMBER = {3},
     PAGES = {223--227},
      ISSN = {0911-0119,1435-5914},
   MRCLASS = {05C35},
  MRNUMBER = {1081196},
       DOI = {10.1007/BF01787573},
       URL = {https://doi.org/10.1007/BF01787573},
}

@article {FF84,
    AUTHOR = {Frankl, P. and F\"uredi, Z.},
     TITLE = {An exact result for {$3$}-graphs},
   JOURNAL = {Discrete Math.},
  FJOURNAL = {Discrete Mathematics},
    VOLUME = {50},
      YEAR = {1984},
    NUMBER = {2-3},
     PAGES = {323--328},
      ISSN = {0012-365X,1872-681X},
   MRCLASS = {05C35 (05C65)},
  MRNUMBER = {753720},
MRREVIEWER = {Ralph\ Faudree},
       DOI = {10.1016/0012-365X(84)90058-X},
       URL = {https://doi.org/10.1016/0012-365X(84)90058-X},
}

@article {KLS18,
    AUTHOR = {Kang, Liying and Liu, Lele and Shan, Erfang},
     TITLE = {Sharp lower bounds for the spectral radius of uniform
              hypergraphs concerning degrees},
   JOURNAL = {Electron. J. Combin.},
  FJOURNAL = {Electronic Journal of Combinatorics},
    VOLUME = {25},
      YEAR = {2018},
    NUMBER = {2},
     PAGES = {Paper No. 2.1, 13},
      ISSN = {1077-8926},
   MRCLASS = {05C50 (05C07 05C65 15A18 15A42)},
  MRNUMBER = {3799419},
       DOI = {10.37236/6644},
       URL = {https://doi.org/10.37236/6644},
}

@incollection {Raz13,
    AUTHOR = {Razborov, Alexander A.},
     TITLE = {Flag algebras: an interim report},
 BOOKTITLE = {The mathematics of {P}aul {E}rd\H os. {II}},
     PAGES = {207--232},
 PUBLISHER = {Springer, New York},
      YEAR = {2013},
      ISBN = {978-1-4614-7253-7; 978-1-4614-7254-4},
   MRCLASS = {05-02 (05C35 05E15)},
  MRNUMBER = {3186665},
MRREVIEWER = {Yandong\ Bai},
       DOI = {10.1007/978-1-4614-7254-4\_16},
       URL = {https://doi.org/10.1007/978-1-4614-7254-4_16},
}

@article {TT22,
    AUTHOR = {Tiner, Gary and Tomlin, Zachery},
     TITLE = {On the {E}rd{\H{o}}s-{S}{\'{o}}s Conjecture for k = 9},
   JOURNAL = {Alabama Journal of Mathematics},
  FJOURNAL = {Alabama Journal of Mathematics},
    VOLUME = {45},
      YEAR = {2022},
     PAGES = {37--45},
       URL = {http://www.ajmonline.org/wp-content/uploads/2022/11/On-the-Erdos-Sos-Conjecture.pdf},
}

@article {Zyk49,
    AUTHOR = {Zykov, A. A.},
     TITLE = {On some properties of linear complexes},
   JOURNAL = {Mat. Sbornik N.S.},
  FJOURNAL = {Mat. Sbornik N.S.},
    VOLUME = {24/66},
      YEAR = {1949},
     PAGES = {163--188},
   MRCLASS = {56.0X},
  MRNUMBER = {35428},
MRREVIEWER = {W.\ T.\ Tutte},
}

@article {Zyk52,
    AUTHOR = {Zykov, A. A.},
     TITLE = {On some properties of linear complexes},
   JOURNAL = {Amer. Math. Soc. Translation},
  FJOURNAL = {Amer. Math. Soc. Translation},
    VOLUME = {1952},
      YEAR = {1952},
    NUMBER = {79},
     PAGES = {33},
   MRCLASS = {56.0X},
  MRNUMBER = {51516},
}

@article {LL81,
    AUTHOR = {Li, Shuo-Yen Robert and Li, Wen Ch'ing Winnie},
     TITLE = {Independence numbers of graphs and generators of ideals},
   JOURNAL = {Combinatorica},
  FJOURNAL = {Combinatorica. An International Journal of the J\'anos Bolyai
              Mathematical Society},
    VOLUME = {1},
      YEAR = {1981},
    NUMBER = {1},
     PAGES = {55--61},
      ISSN = {0209-9683},
   MRCLASS = {05C35},
  MRNUMBER = {602416},
MRREVIEWER = {R.\ L.\ Graham},
       DOI = {10.1007/BF02579177},
       URL = {https://doi.org/10.1007/BF02579177},
}

@incollection {Ste20,
    AUTHOR = {Stein, Maya},
     TITLE = {Tree containment and degree conditions},
 BOOKTITLE = {Discrete mathematics and applications},
    SERIES = {Springer Optim. Appl.},
    VOLUME = {165},
     PAGES = {459--486},
 PUBLISHER = {Springer, Cham},
      YEAR = {[2020] \copyright 2020},
      ISBN = {978-3-030-55857-4; 978-3-030-55856-7},
   MRCLASS = {05C75 (05C05 05C20 05C55 05C65 05C80)},
  MRNUMBER = {4179428},
       DOI = {10.1007/978-3-030-55857-4\_19},
       URL = {https://doi.org/10.1007/978-3-030-55857-4_19},
}

@article {CT91,
    AUTHOR = {Caro, Yair and Tuza, Zsolt},
     TITLE = {Improved lower bounds on {$k$}-independence},
   JOURNAL = {J. Graph Theory},
  FJOURNAL = {Journal of Graph Theory},
    VOLUME = {15},
      YEAR = {1991},
    NUMBER = {1},
     PAGES = {99--107},
      ISSN = {0364-9024,1097-0118},
   MRCLASS = {05C35 (05C65)},
  MRNUMBER = {1090733},
MRREVIEWER = {Michael\ Jacobson},
       DOI = {10.1002/jgt.3190150110},
       URL = {https://doi.org/10.1002/jgt.3190150110},
}

@article {MP07,
    AUTHOR = {Mubayi, Dhruv and Pikhurko, Oleg},
     TITLE = {A new generalization of {M}antel's theorem to {$k$}-graphs},
   JOURNAL = {J. Combin. Theory Ser. B},
  FJOURNAL = {Journal of Combinatorial Theory. Series B},
    VOLUME = {97},
      YEAR = {2007},
    NUMBER = {4},
     PAGES = {669--678},
      ISSN = {0095-8956,1096-0902},
   MRCLASS = {05C35},
  MRNUMBER = {2325805},
MRREVIEWER = {Ivan\ Pashov},
       DOI = {10.1016/j.jctb.2006.11.003},
       URL = {https://doi.org/10.1016/j.jctb.2006.11.003},
}

@article {Nik06,
    AUTHOR = {Nikiforov, Vladimir},
     TITLE = {Walks and the spectral radius of graphs},
   JOURNAL = {Linear Algebra Appl.},
  FJOURNAL = {Linear Algebra and its Applications},
    VOLUME = {418},
      YEAR = {2006},
    NUMBER = {1},
     PAGES = {257--268},
      ISSN = {0024-3795,1873-1856},
   MRCLASS = {05C50},
  MRNUMBER = {2257594},
MRREVIEWER = {Zhongshan\ Li},
       DOI = {10.1016/j.laa.2006.02.003},
       URL = {https://doi.org/10.1016/j.laa.2006.02.003},
}

@article {KN14,
    AUTHOR = {Kang, Liying and Nikiforov, Vladimir},
     TITLE = {Extremal problems for the {$p$}-spectral radius of graphs},
   JOURNAL = {Electron. J. Combin.},
  FJOURNAL = {Electronic Journal of Combinatorics},
    VOLUME = {21},
      YEAR = {2014},
    NUMBER = {3},
     PAGES = {Paper 3.21, 23},
      ISSN = {1077-8926},
   MRCLASS = {05C50 (05C35)},
  MRNUMBER = {3262258},
MRREVIEWER = {Lihua\ Feng},
       DOI = {10.37236/4113},
       URL = {https://doi.org/10.37236/4113},
}

@article {KLM14,
    AUTHOR = {Keevash, Peter and Lenz, John and Mubayi, Dhruv},
     TITLE = {Spectral extremal problems for hypergraphs},
   JOURNAL = {SIAM J. Discrete Math.},
  FJOURNAL = {SIAM Journal on Discrete Mathematics},
    VOLUME = {28},
      YEAR = {2014},
    NUMBER = {4},
     PAGES = {1838--1854},
      ISSN = {0895-4801,1095-7146},
   MRCLASS = {05C50 (05C35)},
  MRNUMBER = {3270186},
MRREVIEWER = {Xiaofeng\ Gu},
       DOI = {10.1137/130929370},
       URL = {https://doi.org/10.1137/130929370},
}

@book {Tur90,
    AUTHOR = {Tur\'an, Paul},
     TITLE = {Collected papers of {P}aul {T}ur\'an. {V}ol. 1--3},
    EDITOR = {Erd\H os, Paul},
 PUBLISHER = {Akad\'emiai Kiad\'o{} (Publishing House of the Hungarian
              Academy of Sciences), Budapest},
      YEAR = {1990},
     PAGES = {Vol. 1: xxxviii+837 pp.; Vol. 2: pp. i--xii and 847--1750;
              Vol. 3: pp. i--xii and 1751--2665},
      ISBN = {963-05-4298-6},
   MRCLASS = {01A75 (11-03)},
  MRNUMBER = {1056293},
MRREVIEWER = {Matti\ Jutila},
}

\appendix
\section{Explicit relation between $r,s$ and $k$ in \cref{thm:Sidorenko-type}}
In this appendix, we will relate positive integers $k,r,s$ with $k\leq r$ satisfying the inequality
\begin{equation}\label{eq:appendix-rks}
    k-s\geq \sum_{i=1}^{s-1}\frac{i}{r-i}.
\end{equation}
We first compute the right hand side.
\begin{lemma}\label{lemma:integral-approx}
    Suppose that $k,r,s$ are positive integers satisfying \cref{eq:appendix-rks}.
    Then $r(k-s) = \Omega(s^2)$, $r-s = \Omega(r)$ and
    \[\sum_{i=1}^{s-1}\frac{i}{r-i} = r\log\left(\frac{r-1}{r-s}\right)-(s-1)+O\left(\frac{s}{r}\right).\]
\end{lemma}
\begin{proof}
    We first show $r(k-s)=\Omega(s^2)$.
    This is clear as 
    \[k-s\geq \sum_{i=1}^{s-1}\frac{i}{r-i}\geq \left\lfloor \frac{s-1}{2}\right\rfloor \frac{\left\lceil \frac{s-1}{2}\right\rceil}{r-\left\lceil \frac{s-1}{2}\right\rceil}=\Omega(s^2r^{-1}).\]
    
    Now we show that $r-s = \Omega(r)$.
    This is clear when $r\geq 2k$, so it suffices to check the case when $r<2k$.
    In this case, we have $2k(k-s)>r(k-s)\geq \Omega(s^2)$.
    This forces $s \leq ck$ for some constant $c<1$, and so $r-s=\Omega(r)$ as $s<k\leq r$.

    Now  let $\cE$ be the error term defined by
    \[\cE =\sum_{i=1}^{s-1}\frac{i}{r-i}-\int_{1}^{s}\frac{x}{r-x}\textup{ d}x = \int_1^s\left(f(\lfloor x\rfloor)-f(x)\right)\textup{ d}x\]
    where we set $f(x) = x/(r-x)$.
    Note that $f'(x) = r(r-x)^{-2}$ is positive and increasing in $x$ when $x\in [1,s]\subseteq [1,r-1]$.
    Therefore
    \[0\geq f(\lfloor x\rfloor)-f(x)\geq (\lfloor x\rfloor -x)f'(x)>-\frac{r}{(r-x)^2}\]
    for any $x\in [1,s]$.
    This shows that 
    \[0\geq \cE \geq -\frac{(s-1)r}{(r-s)^2},\]
    which shows that $\cE = O(s/r)$.
    Therefore
    \[\sum_{i=1}^{s-1}\frac{i}{r-i} =\int_{1}^{s}\frac{x}{r-x}\textup{ d}x+O\left(\frac{s}{r}\right)=  r\log\left(\frac{r-1}{r-s}\right)-(s-1)+O\left(\frac{s}{r}\right).\qedhere\]
    
\end{proof}

\begin{proposition}\label{prop:linear-r}
    Let $r\geq k$ be a positive integer growing with $k$ so that $r = (C+o_{k\to\infty}(1))k$ for some constant $C\geq 1$.
    Then the largest positive integer $s$ satisfying \cref{eq:appendix-rks} also satisfies \linebreak$s=C(1-\exp(-C^{-1})+o_{k\to\infty}(1))k$. 
\end{proposition}
\begin{proof}
    By the choice of $s$, we know
    \[k-s\geq \sum_{i=1}^{s-1}\frac{i}{r-i}\]
    and
    \[k-(s-1)< \sum_{i=1}^{s-2}\frac{i}{r-i}.\]
    Therefore 
    \[k-s+O(1) = \sum_{i=1}^{s-1}\frac{i}{r-i}+O\left(\frac{s}{r-s}\right).\]
    By \cref{lemma:integral-approx}, we know that this implies
    \[k-s +O(1) = r\log\left(\frac{r-1}{r-s}\right)-(s-1)+O\left(\frac{s}{r}\right).\]
    Rearranging, we get
    \[\frac{r-1}{r-s} = \exp\left(\frac{k}{r}+O(r^{-1})\right),\]
    and so
    \begin{align*}
        s =& 1+(r-1)\left(1-\exp\left(-\frac{k}{r}+O(r^{-1})\right)\right)\\
        =& 1+(C+o_{k\to\infty}(1))k \cdot \left(1-\exp\left(-C^{-1}+o_{k\to\infty}(1)\right)\right),
    \end{align*}
    where we use the fact that $r^{-1} = O(k^{-1})$.
    The desired statement thus follows.
\end{proof}

\begin{proposition}\label{prop:quadratic-r}
    Let $d$ be a fixed positive integer.
    Let $k$ be a positive integer with $k>d$, and let $s=k-d$.
    Then the smallest positive integer $r$ satisfying \cref{eq:appendix-rks} also satisfies $r=(\frac{1}{2d}+o_{d;k\to\infty}(1))k^2$.
\end{proposition}
\begin{proof}
    In this proof, we will treat $d$ as a constant and supress all dependencies on $d$.
    
    By the choice of $r$, we also have
    \[d\geq \sum_{i=1}^{s-1}\frac{i}{r-i}\]
    and
    \[d<\sum_{i=1}^{s-1}\frac{i}{r-1-i}.\]
    Note that
    \[\frac{i}{r-i}-\frac{i}{r-1-i} = O(ir^{-2})\]
    for every $i\leq s-1$ as we know that $r-s = \Omega(r)$ by \cref{lemma:integral-approx}.
    Therefore
    \[d = \sum_{i=1}^{s-1}\frac{i}{r-i}+O(s^2r^{-2}).\]
    By \cref{lemma:integral-approx}, we know that $r = \Omega(d^{-1}s^2) = \Omega(s^2)$.
    Therefore  by \cref{lemma:integral-approx},
    \[d = r\log\left(\frac{r-1}{r-s}\right)-(s-1)+O\left(\frac{s}{r}\right)=r\log\left(\frac{r-1}{r-s}\right)-(s-1)+o_{k\to\infty}(1).\]
    
    Note that
    \begin{align*}
        r\log\left(\frac{r-1}{r-s}\right)=&r\log\left(1+\frac{s-1}{r}+\frac{s(s-1)}{r^2}+O\left(\frac{s^3}{r^3}\right)\right)\\
        =&r\left(\frac{s-1}{r}+\frac{(s+1)(s-1)}{2r^2}+O\left(\frac{s^3}{r^3}\right)\right)\\
        =&s-1+\frac{s^2}{2r}+o_{k\to\infty}(1).
    \end{align*}
    Plugging this in, we get
    \[d = \frac{s^2}{2r}+o_{k\to\infty}(1)\]
    and so inverting gives
    \[\frac{2r}{s^2} = \frac{1}{d}+o_{k\to\infty}(1).\]
    Therefore we get that
    \[\frac{r}{k^2} = \left(1+o_{k\to\infty}(1)\right)\frac{r}{s^2}= \frac{1}{2d}+o_{k\to\infty}(1),\]
    as desired.
    
\end{proof}
\end{document}